\documentclass[11pt]{article}

\usepackage{amsmath,amssymb,amsfonts,amsxtra}
\usepackage{verbatim}
\usepackage{diagrams}
\usepackage[all]{xy} 
\usepackage{times}

\diagramstyle[centredisplay,dpi=600,nohug,height=8mm]

\newtheorem{theorem}{Theorem}[section]
\newtheorem{lemma}[theorem]{Lemma}
\newtheorem{corollary}[theorem]{Corollary}
\newtheorem{proposition}[theorem]{Proposition}

\newtheorem{hipo}[theorem]{Hypothesis}

\newtheorem{definition}[theorem]{Definition}
\newtheorem{example}[theorem]{Example}

\newtheorem{paragrafo}[theorem]{}

\newtheorem{remark}[theorem]{Remark}

\numberwithin{equation}{theorem}

\newcommand{\proof}{{\bf{Proof}. }}
\newcommand{\qed}{
$\hfill\square$
}

\newcommand{\BD}{\mathbb D}

\newcommand{\BK}{\mathbb K}

\newcommand{\D}{\mathrm{D}}

\newcommand{\Ch}{\mathsf {Ch}}

\newcommand{\CA}{\mathcal A}
\newcommand{\CB}{\mathcal B}
\newcommand{\CC}{\mathcal C}

\newcommand{\CE}{\mathcal E}
\newcommand{\CF}{\mathcal F}
\newcommand{\CG}{\mathcal G}

\newcommand{\CI}{\mathcal I}

\newcommand{\CL}{\mathcal L}

\newcommand{\CO}{\mathcal O}
\newcommand{\CP}{\mathcal P}

\newcommand{\CR}{\mathcal R}

\newcommand{\lto}{\longrightarrow}

\newcommand{\ot}{\leftarrow}

\newcommand{\curvarrow}[1]{%
\setlength{\unitlength}{0.03\DiagramCellWidth}
\begin{picture}(20,-10)(20,10)
\qbezier (4,26)(25,0)(4,-38)
\put(17.5,-28){\makebox(0,0)[b]{$\scriptstyle {#1}$}}
\put(2.5,-42){\vector(-1,-2){0}}
\end{picture}
}

\DeclareMathOperator{\spec}{Spec}

\DeclareMathOperator{\Hom}{Hom}

\DeclareMathOperator{\id}{id}

\DeclareMathOperator{\ch}{ch}

\DeclareMathOperator{\cone}{C}
\DeclareMathOperator{\cyl}{Cyl}
\DeclareMathOperator{\at}{at}
\DeclareMathOperator{\Lef}{Lef}
\DeclareMathOperator{\At}{At}

\newcommand{\ie}{{\it i.e.} }

\begin{document}

\begin{center}
{\bf\huge{ On the additivity of geometric invariants in derived categories
\\
}}

\ \\

{\bf C. Soneira Calvo}

\end{center}


\hspace{-0,5cm} Departamento de \'Alxebra, Universidade de
Santiago de Compostela,  E-15771 Santiago de Compostela, Spain
(e-mail: carlos.soneira@usc.es)
\ \\


\date{today}

\hyphenation{pseu-do}

\section*{Abstract}

We study the additivity of various geometric invariants involved in Reimann-Roch type formulas and defined via the trace map.
To do so in a general context we prove that given any Grothendieck category $\CA,$ the derived category $\BD(\CA)$ has a compatible triangulation in the sense of [May, J.P. :{\it The Additivity of Traces in Triangulated Categories}, Advances in  Mathematics .{\bf 163}, (2001), 34-73.], but not resorting to model categories, just using the structural properties inherent to $\BD(\CA).$ In the second part of the paper we apply compatibility to prove additivity of traces firstly and then additivity of the Chern character, interpreting this result in terms of a group homomorphism which plays the same role as the Chern character in intersection theory with the $i^{\text{th}}$ Chow group replaced by the $i^{\text{th}}$ Hodge cohomology group.

\bigskip
{\it Keywords: derived category, trace map, Chern character}

{\it 2000 Mathematics Subject Classification: Primary 14C40, 18E15, 18E30, 19L10; Secondary 13D25}

\section*{Introduction}
In order to state Riemann-Roch type formulas in a general algebraic or topological context, as it is proposed, for example in \cite{FM}, additivity of different invariants is needed, in particular Chern and Leftschetz characters (the last one for bivariant versions of the theorem), defined via the trace map. Since in both algebraic and topological cases,  derived categories arise naturally, we study this case, where these invariants are defined via generalized trace maps.This being the situation, the main problem is that in  the derived category, the trace map is not additive (see \cite{Fe} for an easy example).
In the enlightening paper \cite{May}, are explained five axioms, that we also refer as (TC1) - (TC5), in the triangulated structure of the category, guarantying additivity of traces. This is what May calls a compatible triangulation.
A generic case covering our expectative of generality is that of $\BD(\CA)$ with $\CA$ a Grothendieck category, because it embraces the general framework of algebraic geometry.
In the first part of this work it is proved the theorem below, being the technical tool we will apply in the remaining to deal with geometric invariants.

\begin{theorem}\label{th}
The category $\BD(\CA)$ has a triangulation compatible with its closed monoidal structure, and it can be given explicitly.
\end{theorem}

In the second part this is applied to additivity of invariants.

In \cite{May} both proofs are done for a category $\CC$ with model category $\CB$ such that the closed symmetric monoidal structure on $\CC$ is induced from one in $\CB.$
We give a direct proof of Theorem \ref{th} using the very structural properties of the derived category, without resorting to model category theory, reducing the question to a suitable and, to a certain extent, natural,  explicit  choice of morphisms in the category of chain complexes.
Speaking figuratively, we would say that  the category of chain complexes becames our model.
Then this new structure provided by the axioms is used to prove additivity.


In the first section we fix the general hypothesis of this work, giving some familiar examples where they are fullfilled, pointing out that they can be assumed in most situations arising in algebraic geometry. We recall also some technical results about abstract duality that will be used later.

In the second one, we recall the notion of compatible triangulation; axioms (TC1), (TC2), (TC3) and (TC4) are stated and proved, while axiom (TC5) and the the main result of the first part of the paper, Theorem \ref{th}, are just enunciated.  We must remark that techniques used in this work are essentially different from those appearing in \cite{May}, because we do not assume the existence of a model category; our tools are homological and homotopical algebra, focusing on the triangulated structure of the derived category, heavily using the different implications of the octahedron axiom, the  properties of derived functors and those of  involution.

In the third section we prove axiom (TC5). We do it in a separated section due to its length and because the techniques required are different. It rests on two technical lemmas about abstract duality and we argue by giving  explicitly
the morphisms whose existence is claimed in (TC5). In the proof of axiom (TC5) we also use the previous axioms, being relevant to point out that, although those were already proved it is not enough to make this procedure logically consistent, it is also needed that the successive choices of morphisms we must do in each axiom to be compatible to the preceding ones. Here the key point is that we always take morphisms arising from the inherent triangulated structure of the category such as completion to an octahedron.

After that, it begins the second part of the paper, applying Theorem \ref{th} to prove the additivity of geometric invariants, that are based on the generalized trace maps.

In the fourth section we briefly define the trace map viewing it as a particular case of an orientation and prove its additivity in Theorem \ref{addtraces}. In our hypothesis, we can't apply the homotopy extension property, which is the main ingredient in the proof of additivity in model category framework. The main  problem is  to find an arrow with certain properties making it natural in some sense. In $\BD(\CA)$ or $\BK(\CA)$, we can always assume the existence of certain arrows but we there is not naturality.
Our strategy consists on prove Lemma \ref{lemacyl} before and use it as a tool, in the sense that ir allows us to choose morphisms coming from arrows in the category of chain complexes and avoid any compatibility problem. Intuitively, what we do in Lemma \ref{lemacyl} is to use the cylinder construction to codify the homotopy data and work in $\Ch(\CA)$ where morphisms can be given in a natural way. This procedure makes also possible to give the required morphisms explicitly.

The last aim is to prove the additivity os Chern character, so in Section 5 we restrict to the derived category of quasi-coherent sheaves of modules over a scheme $X$. After reviewing the definition and basic functorial properties of the Atiyah class and higher Atiyah classes, we apply Theorem \ref{addtraces} and the fact that two isomorphic objects in $\BD(\CA_{qc}(X))$ define the same Chern character in each degree, proved in Proposition \ref{invqischi}, to demonstrate the following

\begin{theorem}\label{addchi} Let
\[\CE\stackrel{f}{\lto}\CF\stackrel{g}{\lto}\CG\stackrel{h}{\lto}\Sigma\CE\]
be a distinguished triangle. Then
\[\ch_i(\CF)=\ch_i(\CE)+\ch_i(\CG)\]
for every $i\in\mathbb{N}$.
\end{theorem}

It is now easy to obtain, via the definition of the total Atiyah class, an analogous result of additivity for a Chern character with values in Hodge cohomology defined as
\[\ch(\CE):=\Lef(\id_\CE,\At):\CO_X\lto\bigoplus_{i=0}^n\Sigma^i\Omega_{X|S}^i.\]
 where $\ch(\CE)=\sum_{i=0}^n\frac{1}{i!}\ch_i(\CE).$ Here $n$ is the dimension of the scheme over the base object of the category and we assume that $n!$ is an unity in the sheaf of sections of the scheme base.

Finally we organize the preceding results as a map playing a central role in Riemman-Roch theorems as claimed at the beginning of the present introduction; particulary there is a homomorphism of groups

\[\ch:K^0(X)\lto\bigoplus_{i=0}^n H^i(X,\Omega_{X|S}^i).\]

\section{Notations and conventions}

 From now on let $\CA$ be a  Grothendieck category.
 We will consider its associated derived category  $\BD(\CA).$ By the results in \cite{Lip}  it has $q$-injective reolution, therefore functors on the homotopy categoy $\BK(\CA)$ induced by left exact functors from $\CA$ itself admit a right derived functor. Also, Brown representability holds. We will denote the suspension functor in $\BK(\CA)$ and $\BD(\CA)$ by $\Sigma,$ while by $\Sigma^{-1}$ we denote the desuspension functor. The category of complexes of objects of $\CA$ will be denoted by $\mathbf{C}(\CA).$ We fix once and for all that if $E\in\mathbf{C}(\CA),$ its differential will be denoted by $d_E.$

 For a map $f: E\lto F$ has associated a complex, its cone $\cone(f)$, that fits into an exact sequence
 \[0\lto F\lto\cone(f)\lto\Sigma E\lto0\]
in $\mathbf{C}(\CA)$ that yields a distinguised triangle
\[E\stackrel{f}{\lto}F\lto\cone(f)\lto \Sigma E\]
in $\BK(\CA)$ and $\BD(\CA).$ See \cite[Section 1.3 and Examples 1.4.4]{Lip} for the relevant sign conventions.

We will consider in addition that it has a closed structure in the sense of \cite{Ei}. It means that it is symmetric monoidal with an internal tensor bifunctor that we will denote by $-\otimes-$. In addition, it has an internal hom functor left adjoint to the tensor (with the first and the second variable fixed for tensor and internal hom, respectively). We will denote the internal hom by $\mathcal{H}om(-,-).$ Note that this makes $\otimes$ a right exact functor (in either variable by symmetry) and $\mathcal{H}om(-,-)$ an exact functor in the second variable. The base object will be denoted by $S.$

Both functors extend to the homotopy category. The derived functor of $\mathcal{H}om(-,-)$ will be denoted as $\mathbf{R}\mathcal{H}om(-,-)$
 but most of the time we will abbreviate it by $[-,-].$

 It may not be the case that the functor $\otimes$ is exact so we make the following hypothesis.

\begin{hipo}\label{hipo}
The category $\CA$ possesses $q$-flat resolutions in the sense of \cite[Definition 2.5.1]{Lip} .
\end{hipo}

Under this hypothesis the derived functor of $-\otimes-$ exists and we will denote it by $-\otimes-$ by simplicity. It makes $\BD(\CA)$ a symmetric monoidal category with unit object $S$ seen as a complex concentrated in degree $0$. The functor $[-,-]$ is the remaining ingredient that makes  $\BD(\CA)$ a closed category.

\begin{example}\label{exemplos} Our assumptions are natural in the sense that they are  widely verified in the general framework of algebraic geometry.

 \begin{enumerate}

 \item[(i)] The category $\CA=R-Mod$ with $R$ a commutative ring.

 \item[(ii)] The category $\CA=\CR-Mod$  $(X,\CR)$ a ringed space. The fact that $\CR-Mod$ is Grothendieck is a classical fact, for the existence of $q$-flat resolutions, see  \cite{Sp}.

 \item[(iii)]  The category $\CA=\CA_{qc}(X)$ of sheaves of quasi-coherent $\CO_X$-Modules with  $(X,\CO_X)$ a quasi compact and separated scheme, see \cite{Se}.

\item[(iv)] Let $B$ be a commutative Hopf Algebra, $\CA=B-Comod$ the category of left $B$-comodules. It falls under our conditions essentially by \cite[\S 9.5]{HPS}.

\end{enumerate}

\end{example}

\begin{remark}
 Note that examples (ii), (iii) and (iv) with a twist are different generalizations of example (i). For (iv), look at \cite[Lemma 9.3.5(b)]{HPS}.

\end{remark}

\begin{paragrafo}\label{absdual}\bf{Review of abstract duality}
\end{paragrafo}
We summarize now the elements of duality in the context of a closed category. We follow the treatment of \cite{LMS} but using the conventions of $\cite{Lip}$ especially for the internal composition, and giving explicitly the maps involved and refer to \cite{LMS} for most proofs.

Given objects $E,F\in\BD(\CA)$, consider the canonical internal evaluation map

\[e_{E,F}:[E,S]\otimes F\lto [E,F]\]
induced by
\[e_{E,S,F}:[E,S]\otimes [S,F]\lto[E,F]\]
from \cite[Exercise (3.5.3)(c)]{Lip}  through the natural map
\[F\lto[S,F]\]
adjoint to \cite[Definition (3.4.1)]{Lip}
\[\rho:F\otimes S\lto F.\]

The object $E\in\BD(\CA)$ is said to be strongly dualizable if, and only if, $e_{E,F}$ is an isomorphism for all $F\in\BD(\CA)$.

We will use the notation $\mathrm{D}E:=[E,S]$ for duality of strongly dualizable objects.

There is a canonical morphism
\[t_E:\mathrm{D}E\otimes E\lto S\]
as in \cite[Exercise (3.5.3)(b)]{Lip}  that induces
\[\beta:E\lto\mathrm{D}\mathrm{D}E.\]

\begin{lemma}\label{ddual} The morphism $\beta$ is an isomorphism whenever $E$ is strongly dualizable.
\end{lemma}
\proof See \cite[Proposition III.1.3(i)]{LMS}.
\qed

Denote by
\[u_E:S\lto E\otimes\mathrm{D}E\]
the canonical morphism induced from the one in \cite[Exercise (3.5.3)(f)]{Lip}.
We have the following useful fact

\begin{lemma}\label{diagdual} Let $f:E\lto F$ a map in $\BD(\CA)$ with $E$ and $F$ strongly dualizable. The diagrams
\begin{diagram}
\begin{diagram}
\mathrm{D}F\otimes E & \rTo^{\id\otimes f} & \mathrm{D}F\otimes F \\
\dTo^{\mathrm{D}f\otimes\id} &    & \dTo_{t_F}  \\
\mathrm{D}E\otimes E & \rTo_{t_E} & S
\end{diagram}
&
\begin{diagram}
 S & \rTo^{u_F}&F\otimes\mathrm{D}F \\
\dTo^{u_E} &    & \dTo_{id\ot\mathrm{D}f}  \\
E\otimes\mathrm{D}E & \rTo_{f\otimes id} & \mathrm{D}F\otimes E
\end{diagram}
\end{diagram}
commute.
\end{lemma}
\proof See \cite[Proposition III.1.5]{LMS}. \qed

\section{May's Axioms}

When a closed category has an additional structure the problem of expressiong the compatibility of the closed structure with the additional one arises. In the case of an abelian category, the exactness properties of the tensor and internal hom provide all the compatibilities with short exact sequences one could need.

The case of a triangulated category is subtler because the triangulation is an additional structure. On the contrary, in the case of abelian categories, exact sequences arise from the properties satisfied by the category.

May proposed in \cite{May} a series of axioms that try to express this compatibility. In the paper the author showed how to check them under the hypothesi that the triangulated category $\mathbf{T}$ possesses a model in the sense of Quillen. In our case a Quillen model for $\BD(\CA)$ is known to exist but the interplay between this model and the closed structure seems difficult to use, therefore we approach the proof in an independent fashion, exploiting the underlying additive structure of $\BD(\CA)$.

Let us state May's axioms and show how to prove the straightforward ones. For objects $E,F\in\BD(\CA)$ the symmetry or switch map will be denoted by
\[\gamma=\gamma_{E,F}:E\otimes F\lto F\otimes E\]
See \cite[Definition 3.4.1]{Lip}. Most of the time we will omit the subscript as it should be clear from the context.

{\bf Axiom (TC1)}  Let $E\in\BD(\CA)$. There is a natural isomorphism $\alpha: E\otimes  \Sigma S \lto \Sigma E$ such that the composite

\[\Sigma\Sigma S \stackrel{\alpha^{-1}}{\lto}\Sigma S\otimes \Sigma S\stackrel{\gamma}{\lto} \Sigma S\otimes \Sigma S \stackrel{\alpha}{\lto}\Sigma\Sigma S\]

is multiplication by $-1$, where $\gamma$ is the symmetry map \cite[Definition 3.4.1]{Lip}.

{\bf Axiom (TC2)} The functors $[-,-]$ and $-\otimes-$ are $\Delta$-functors \cite{Lip} \S1.5 in both variables.

\begin{remark}
The content of the previous axioms is essentially equivalent to a compatible closed structure on a triangulated category in the sense of Hovey, Palmieri and Strickland \cite[Definition A.2.1]{HPS}.
\end{remark}


{\bf Axiom (TC3)} Given distinguished triangles
\[E  \stackrel{f}{\lto} F  \stackrel{g}{\lto} G  \stackrel{h}{\lto}\Sigma E\]
and
\[E'  \stackrel{f'}{\lto} F'  \stackrel{g'}{\lto} G'  \stackrel{h'}{\lto}\Sigma E',\]
 there is an object $W$ and morphisms $p_i$ and $j_i$ with $i=1,2,3;$ such that the following triangles are distinguished
\[F\otimes E'  \stackrel{p_1}{\lto} W  \stackrel{j_1}{\lto} E\otimes G'  \stackrel{f\otimes h'}{\lto}\Sigma (F\otimes E')\]

\[\Sigma^{-1}(G\otimes G')  \stackrel{p_2}{\lto} W  \stackrel{j_2}{\lto} F\otimes F'  \xrightarrow{-g\otimes g'}
G\otimes G'\]

\[E\otimes F'  \stackrel{p_3}{\lto} W  \stackrel{j_3}{\lto} G\otimes E'  \stackrel{h\otimes f'}{\lto}\Sigma (E\otimes F')\]
and $W$ is the common apex of three octahedron arising by completing  triangles

 \begin{tabular}{ccc}
\xymatrix@=6mm@W=0.5mm{\Sigma^{-1}(G\otimes G')\ar[rr]\ar[rd]_{p_2} & &F\otimes F'\\
& W\ar[ru]_{j_2}&
} &\xymatrix@=6mm@W=0.5mm{F\otimes E'\ar[rr]\ar[rd]_{p_1} & &E\otimes G'\\
& W\ar[ru]_{j_1}&
}\end{tabular}

braided in the following way

{\scriptsize\begin{diagram}
E\otimes G'     &        &\lTo^{h\otimes\id} &        &\Sigma^{-1}G\otimes G'       &          & \rTo^{\Sigma^{-1}h\otimes\id}  &         &  E\otimes G'   & \lTo^{\id\otimes g'}   &     &         &E\otimes F' \\
&\luTo  &  \circlearrowleft &\ldTo & & \rdTo &  \circlearrowleft &\ruTo &              & \luTo& \circlearrowleft &\ldTo &           \\
\dTo^{f\otimes h'}    & +     & W  &   +      &  \uTo^{-g\otimes g'}    & +          &   W  & +         &  \dTo^{f\otimes h'}  & +        & W &    +    & \uTo_{h\otimes f'}           \\
       &\ruTo& \circlearrowleft &\rdTo &         & \ldTo  & \circlearrowleft &  \luTo&                        &\ruTo & \circlearrowleft  & \rdTo&           \\
F\otimes E'   &        &\rTo_{\id\otimes f'} &        &F\otimes F'       &         &\lTo_{\id\otimes f'}   &          &F\otimes E'    & \rTo_{g \otimes\id}   &     &         &  G\otimes E'
\end{diagram}}

where the  triangles marked $\circlearrowleft$ are commutative and the ones marked $+$ distinguished. Note that this is equivalent to May's formulation, {\it cf.} \cite[page 49]{May}.

In \cite{Kell-Nee},  the tensor product of the ambient category is required to be  decent, and it is an apparently new concept. Although in our context it is easy to prove that tensor product verifies similar properties to a decent one,  we focus on the triangulated structure of $\BD(\CA),$ coming from the one on  $\BK(\CA).$ In this way, and differing from \cite{Kell-Nee}, unbounded complexes are allowed.

\begin{lemma}\label{(TC3')} Consider again the distinguised triangles
\[E  \stackrel{f}{\lto} F  \stackrel{g}{\lto} G  \stackrel{h}{\lto}\Sigma E\]
and
\[E'  \stackrel{f'}{\lto} F'  \stackrel{g'}{\lto} G'  \stackrel{h'}{\lto}\Sigma E',\]
 there is an object $V$ and morphisms $k_i$ and $q_i$ with $i=1,2,3;$ such that the following triangles are distinguished
\[E\otimes G^\prime  \stackrel{k_1}{\lto} V  \stackrel{q_1}{\lto} G\otimes F^\prime  \stackrel{f\otimes h^\prime}{\lto}\Sigma (E\otimes G^\prime)\]

\[\Sigma^{-1}(F\otimes F^\prime)  \stackrel{k_2}{\lto} V  \stackrel{q_2}{\lto} E\otimes E^\prime  \xrightarrow{-g\otimes g^\prime}
F\otimes F^\prime\]

\[G\otimes E^\prime  \stackrel{k_3}{\lto} V  \stackrel{q_3}{\lto} F\otimes G^\prime  \stackrel{h\otimes f^\prime}{\lto}\Sigma (G\otimes E^\prime)\]
and $V$ is the common apex of three octahedron arising by completing  the commutative triangles, analogously to the case of (TC3).
\end{lemma}

\proof Follows immediately from axiom (TC3), see \cite{May}, Lemma 4.7.\qed

\begin{paragrafo}\label{barra}{\bf Involution process}.
Reversing the order of the distinguished triangles $(f,g,h)$ and $(f',g',h')$ and applying (TC3) and (TC3') we obtain corresponding objects $\overline{W}$ and $\overline{V}.$ Also, by the naturality of the switch transformation $\gamma$ and completing the triangles, we get isomorphisms $\overline{\gamma}:W\rightarrow\overline{W}$ and $\overline{\gamma}:V\rightarrow\overline{V}$ such that
\[\overline{\gamma}p_2=\overline{p_2}\Sigma^{-1}\gamma,\qquad \overline{j_2}\gamma=\overline{\gamma}j_2,\qquad\overline{\gamma}k_2=\overline{k_2}\Sigma^{-1}\gamma,\qquad \overline{q_2}\gamma=\overline{\gamma}q_2.\]
Then, for every $r\in\{p_i, j_i, k_i, q_i\text{ with }i=1,2,3\}$ we define the corresponding $\overline{r}:=\overline{\gamma}r\gamma^{-1}$ with a negative shift if necessary. With these definitions, there are commutative diagrams as in Axioms (TC3) and (TC3') for the interchanged triangles. We call this operation an involution of the original diagrams. See \cite[Remark 4.10]{May}.
\end{paragrafo}

{\bf Axiom (TC4)} (\emph{The additivity Axiom}). If $W$ is the common apex arising from (TC3) and $V$ the one from (TC3'), and  maps $j_i$ and $k_i$ as brefore above, the
\[ V\stackrel{(j_2,(j_1,j_3))}{\lto} (F\otimes F')\oplus(E\otimes G')\oplus(G\otimes E') \stackrel{(k_2,(k_1,k_3))}{\lto} W\lto \Sigma V  \]
is a distinguished triangle.

 {\bf Axiom (TC5)}\label{TC5} (\emph{The Braid Duality Axiom}). There is a choice for the morphism $k_i$ with $i=1,2,3$, arising when apply the construction of Axiom (TC3') to the distinguished triangles
\[\mathrm{D}G\stackrel{\mathrm{D}g}{\lto}\mathrm{D}F\stackrel{\mathrm{D}f}{\lto}\mathrm{D}E\stackrel{D\Sigma^{-1}h}{\lto}\Sigma \mathrm{D}G\]
and
\[E\stackrel{f}{\lto}F\stackrel{g}{\lto}G\stackrel{h}{\lto}\Sigma E\]
such that:

(a)  There is a map $\overline{t}:\overline{W}\lto S$ making the diagram
\begin{diagram}
(\mathrm{D}G\otimes G)\oplus(\mathrm{D}E\otimes E) & \rTo^{\overline{k_1},\overline{k_3}} &\overline{W}&\lTo^{\overline{k_2}}& \mathrm{D}F\otimes F\\
 & \rdTo_{(t_G,t_E)} &\dTo\overline{t}&\ldTo_{t_F}& \\
 &                                          &    S                 &                        &
\end{diagram}
commutative. Where $t$ denote the evaluation maps as in \ref{absdual}.

(b) If $E$, $F$ and $G$ are strongly dualizable, then, with a choice of morphisms for the diagram of (TC3') with respect to the   indicated triangles, make this (TC3') type diagram isomorphic to the (TC3) type diagram resulting make this (TC3') diagram isomorphic to the (TC3) diagram that arises when applying axiom (TC3) to the triangles:
\[E\stackrel{f}{\lto} F\stackrel{g}{\lto} G\stackrel{h}{\lto} \Sigma E\]
and
\[\mathrm{D}G\stackrel{\mathrm{D}g}{\lto}\mathrm{D}F\stackrel{\mathrm{D}f}{\lto}\mathrm{D}E\stackrel{\mathrm{D}\Sigma q^{-1}h}{\lto}\Sigma \mathrm{D}G,\]
and Axiom (TC4) is satisfied with respect to the above (TC3') type diagram and an involution of the latter (TC3) type diagram.

\begin{definition}
The closed monoidal structure of $\BD(\CA)$ is said to be compatible (with its canonical triangulated structure) if the five Axioms (TC1)-(TC5), just stated, are satisfied.
\end{definition}

\begin{lemma}\label{tc1h}
Axiom (TC1) holds.
\end{lemma}

\proof
The object $S$  is  the complex being is $S$ of $\CA$ in degree zero and zero in any other degree.
Then, if $E\in\BD(\CA)$,
\[{(E\otimes \Sigma S)}^n=\oplus_{i+j=n}(E^i\otimes S^j)=E^{n+1}\]

because $S^j=0$ unless $j=0,$ and the differential is
\[\sum_{i+j=n} ({(-1)}^i{\delta_A}^i\otimes {\delta_{\Sigma S}}^j ) + ({\delta_A})^j\otimes \id={\delta_A}^{n+1}\otimes\id\]
and we can define $\alpha$ as the morphism of complexes that in each degree is ${(-1)}^n\cdot\id$, so it happens that $\alpha=\alpha^{-1}$.
The symmetry morphism $\gamma: \Sigma S\otimes \Sigma S\lto \Sigma S\otimes \Sigma S$ is multiplication by $-1$, because when we calculate the total complex, the vertical differentials become horizontal differentials.

Then, putting together these two descriptions, we conclude that the isomorphism of Axiom (TC1) is multiplication by $-1$, as claimed.
\qed

\begin{lemma}\label{tc2h}
Axiom (TC2) holds.
\end{lemma}

\proof
Trivial, both are derived functors of $\Delta$-bifunctors in $\BK(\CA).$ \qed

\begin{proposition}\label{tc3h}
Axiom (TC3) holds.
\end{proposition}

\proof

We represent both distinguished triangles by semi-split  exact sequences in $\BK(\CA),$namely
\[0\lto E \lto F \lto G \lto 0\]
\[0\lto E' \lto F' \lto G' \lto 0\]
Taking the tensor product of both sequences we obtain a $3\times 3$ diagram with exact and semi-split rows and columns

{\begin{diagram}
   &     & 0               &      &  0              &     & 0               &        &     \\
   &     &\dTo           &      &\dTo            &    &\dTo            &        &      \\
0 &\rTo&E\otimes E'&\rTo&F\otimes E'&\rTo&G\otimes E'& \rTo  &0\\
   &     &\dTo           &      &\dTo            &    &\dTo            &       &      \\
0 &\rTo&E\otimes F'&\rTo&F\otimes F'&\rTo&G\otimes F' &\rTo  &0\\
  &     &\dTo           &      &\dTo            &    &\dTo            &        &      \\
0 &\rTo&E\otimes G'&\rTo&F\otimes G'&\rTo&G\otimes G' &\rTo  &0\\
  &     &\dTo           &      &\dTo            &    &\dTo            &        &      \\
 &     & 0               &      &  0              &     & 0               &        &     \\
\end{diagram}}

From this set-up we argue much as the same as in \cite{Kell-Nee}. However on this paper the tensor product of the ambient category is required to be  decent, and it is an apparently new concept. Although in our context it is easy to prove that tensor product verifies similar properties to a decent one,  we focus on the triangulated structure of $\BD(\CA),$ coming from the one on  $\BK(\CA).$ They also insist of the category that models the triangles being bounded, but we deal with unbounded complexes as well. For the reader's convenience we sketch how to paste the present setup with the one at {\it loc.cit.}

Define $V:=(F\otimes F')/(E \otimes E')$ as an object of $\mathbf{C}(\CA).$ We have a distinguished triangle
\[E\otimes E'\lto{} F\otimes F'\lto{} V\stackrel{+}{\lto}\]
From here, by a diagram chase argument we deduce an exact diagram (in $\mathrm{C}(\CA)$)
\begin{diagram}
   &     & 0                &      &  0              &      & 0               &        &     \\
   &     &\dTo            &      &\dTo            &     &\dTo            &        &      \\
0 &\rTo&E\otimes E'&\rTo&F\otimes E'&\rTo&G\otimes E' &\rTo  &0\\
   &     &\dTo_{\id}    &      &\dTo            &      &\dTo            &       &      \\
0 &\rTo&E\otimes E'&\rTo&F\otimes F'&\rTo        &  V            & \rTo  &0\\
  &     &\dTo             &      &\dTo            &              &\dTo            &        &      \\
 &      & 0                 &\rTo&F\otimes G'&\rTo_\id   &F\otimes G' &\rTo  &0\\
  &     &                   &      &\dTo            &    &\dTo            &        &      \\
 &     &                    &      &  0              &     & 0               &        &     \\
\end{diagram}

And a semi-split exact sequence

\[0\lto G\otimes E' \lto V \lto F\otimes G' \lto 0.\]

Similarly, we obtain

\[0\lto E\otimes G' \lto V \lto G\otimes F' \lto 0.\]

yielding distinguished triangles
\[ \Sigma^{-1}F\otimes G' \lto G\otimes E'\lto V\lto F\otimes G' \]
and
\[\Sigma^{-1}G\otimes F'\lto E\otimes  G'\lto V\lto G\otimes F'.\]

Summing up,  $V$ is the common cone of the three diagonal arrows on the following diagram in $\BD(\CA)$
\begin{diagram}
 &   &                    &       &\Sigma^{-1}(F\otimes G') &\rTo & \Sigma^{-1}(G\otimes G')\\
 &   &                    &       &    \dTo                              &\rdTo&\dTo  \\
 &   & E\otimes E' &\rTo & F\otimes E'                      &\rTo  &G\otimes E' \\
 &   & \dTo            &\rdTo&\dTo                                 &        &\dTo  \\
\Sigma^{-1} G\otimes F'&\rTo   & E\otimes F' &\rTo & F\otimes F'                      &\rTo  &G\otimes F' \\
\dTo& \rdTo& \dTo &  & \dTo &  &\dTo \\
\Sigma^{-1} G\otimes G'&\rTo   & E\otimes G' &\rTo & F\otimes G'&
\rTo  &G\otimes G'
\end{diagram}

Consider now the two octahedra obtained by completing the triangles occurring in the commutative square from which $V$ was defined. In the octahedron, both squares containing the apices can be chosen to be homotopy pull-backs (see \cite{Nee2} for a discussion of this concept), so the object $V$ and arrows  involved constitute a (TC3') type diagram. By symmetry, we know that if we do the construction beginning from the original triangles applying first the suspension functor, we will get an object $W$ and the corresponding morphisms that make axiom (TC3) hold. See \cite[Pages 541-547]{Kell-Nee} for further details. \qed

\begin{lemma}\label{tc4h}
Axiom (TC4) holds.
\end{lemma}

\proof
By \cite[Theorem 4.1]{Kell-Nee} the axiom (TC4) is a formal consequence ofAxiom (TC3). \qed

{\it Remark} In the same proof they show that

\[W\lto(F\otimes E')\oplus(E\otimes G')\lto V\lto \Sigma W\]

is a distinguished triangle.

We delay the proof of axiom (TC5) to the next paragraph, however we state now our first main theorem.

\begin{theorem}\label{th_ct}
Under the general assumptions and assuming hypothesis \ref{hipo} holds, the closed monoidal structure of $\BD(\CA)$ is compatible.The compatibility maps can be chosen explicitly.
\end{theorem}

\proof
Combine lemmas (\ref{tc1h}), (\ref{tc2h}) and (\ref{tc4h}) whit Proposition (\ref{tc4h}) and Theorem (\ref{tc5h}) below.  \qed

\section{Strongly dualizable objects and braid duality}

This section is devoted to the verification of Axiom (TC5). We use a separate section due to its length and the different flavor of the techniques employed. As the objects involved are strongly dualizable we will use frequently the results reviewed in \ref{absdual}. We begin with a couple of Lemmas that depend on the strongly dualizable hypothesis and that will be often invoked throughout the proof.

In he previous setting, let $E,G\in\BD(\CA)$ be strongly dualizable objects. There is a canonical map
\[p:\mathrm{E}\otimes\mathrm{G}\lto\mathrm{E\otimes G}\]
composing \cite{Lip} Exercise (3.5.3)(c) with the canonical isomorphism $[S,S]\cong S.$ Let $E,F\in\BD(\CA)$ be strongly dualizable objects. We define the morphism $\xi$ as the composition
\[\mathrm{E}\otimes F\stackrel{id\otimes\beta}{\lto}\mathrm{D}E\otimes \mathrm{D}\mathrm{D}F\stackrel{p}{\lto}\mathrm{D}(E\otimes\mathrm{D}F)\]
Where in the morphism $p,$ $G=\mathrm{D}F$ and $\beta:F\stackrel{\sim}{\lto}\mathrm{D}\mathrm{D}F$ is the isomorphism from Lemma \ref{ddual}.

\begin{lemma}
The morphism $\xi$ is an isomorphism.
\end{lemma}

\proof it is clear that $id\otimes\beta$ is an isomorphism and so is $p$ by \cite{LMS}, Proposition III.1.3(iii), beign $F$ and $G$ strongly dualizable.
\qed.

\begin{lemma}\label{funxi} Consider the isomorphism $\xi$ and let $f:E\rightarrow E'$ and $g:F\rightarrow F'$ be morphisms between strongly dualizable objects, then

\begin{itemize}
\item[(i)] $\mathrm{D}(f\otimes id_{\mathrm{D}F})\xi=\xi(\mathrm{D}F\otimes id_F)$

\item[(ii)]$\mathrm{D}(id_E\otimes\mathrm{D}g)\xi=\xi(id_E\otimes g)$
\end{itemize}

\end{lemma}
\proof For (i) the functoriality of $\beta$ yields the following commutative diagram:

\[\begin{diagram}
\mathrm{D}E'\otimes F & \rTo^{\id\otimes\beta} & \mathrm{D}E'\otimes \mathrm{D}\mathrm{D}F &\rTo^p &     \mathrm{D}(E'\otimes \mathrm{D}F)\\
\dTo_{\mathrm{D}f\otimes\id} &    & \dTo_{\mathrm{D}f\otimes\id}&  &\dTo_{\mathrm{D}(f\otimes\id)}\\
\mathrm{D}E\otimes F & \rTo_{\id\otimes\beta} & \mathrm{D}A\otimes \mathrm{D}\mathrm{D}F &\rTo_p &         \mathrm{D}(E\otimes \mathrm{D}F)
\end{diagram}\]

Part (ii) is similar.
\qed

\begin{theorem}\label{tc5h}
Axiom (TC5) holds.
\end{theorem}

\proof
We consider firs the proof of part (a). Take a semi-split sequence in $\mathbf{C}(\CA)$

\begin{equation}\label{EFG}
0\lto E\stackrel{f}{\lto}F\stackrel{g}{\lto}G\lto0
\end{equation}

Consider the functor $\mathrm{D}=[-,S]$ described by $Hom(-,I_S)$ with $I_S$ a $q$-injective resolution of the unit object $S.$ The sequence
\[0\lto \mathrm{D}G\stackrel{\mathrm{D}g}{\lto}\mathrm{D}F\stackrel{\mathrm{D}f}{\lto}\mathrm{D}E\lto0\]
is again semi-split, and therefore the triangles
\[ E\stackrel{f}{\lto}F\stackrel{g}{\lto}G\stackrel{h}{\lto}\Sigma E\]

\[ \mathrm{D}G\stackrel{\mathrm{D}g}{\lto}\mathrm{D}F\stackrel{\mathrm{D}f}{\lto}\mathrm{D}E\stackrel{D\Sigma^{-1}h}{\lto}\Sigma \mathrm{D}G\]

are distinguished, where $h$ denotes the homotopy invariant of the semi-split sequence (\ref{EFG}).

Take the distinguished triangle
\[\D G\otimes E\stackrel{\D g\otimes f}{\lto}\D F\otimes F\lto W\stackrel{+}{\lto}\]
\ie $W=:=\cone(\D g\otimes f)$
and consider the diagram \label{diagram (TC5)(a)}
\begin{diagram}
                   &                       &  \D G\otimes E                           &                      &\\
                   &\ldTo^{\id\otimes f}&     &\rdTo^{\D g\otimes\id}   &   \\
\D G\otimes F&\rTo^{\D g\otimes\id}&\D F\otimes F&\lTo^{\id\otimes f}&\D F\otimes E\\
\dTo^{\id\times g}&   &\dDashto^{k_2}  & \curvarrow{t_F} & \dTo_{\D f\otimes\id}\\
\D G\otimes G&\rDashto^{k_1}& W &\lDashto^{k_3}         &\D E\otimes E\\
                   &\rdTo_{t_G}&\dDashto^{\overline{t}}&\ldTo_{t_E}&   \\
                   &                       &  S                           &                      &
\end{diagram}

We will show that there are morphisms $k_1,$ $ k_2,$ $ k_3$ and $\overline{t}$ making the whole diagram commutative. The roof is clearly commutative. Also the two squares formed with the curved arrow are commutative by Lemma \ref{diagdual}.

Our next task is to construct (and fix) a morphism $\overline{t}.$ Note that the composition
$t_F(\D g\otimes f) : \D G\otimes E\lto S$ is the zero morphism, because
 \[t_F(\D g\otimes f)=t_F(\D g\otimes id)(id\otimes f)=t_G(id\otimes g)(id\otimes f)= t_G(id\otimes gf)=0,\]
 where in the second equality we used Lemma \ref{diagdual}.

As in Axiom (TC3) take $W=\cone(\D g\ot f)$ and consider the diagram

\begin{diagram}
\D G\otimes E&\rTo^{} & \D F\otimes F &\rTo^{k_2} & W &\rTo&\Sigma(\D G\otimes E)\\
\dTo^{}  &        & \dTo^{t_F}  &              &\dDashto_{\overline{t}}      &        &\dTo_{}\\
0&\rTo & S&\rTo_{\id} &S & \rTo&0
\end{diagram}
Note that, as a matrix, $\overline{t}=(0\quad t_F)$ and $\overline{t}k_2=t_F.$

Consider now the morphism of triangles

\begin{diagram}
\D G\otimes E&\rTo^{\id\otimes f} & \D G\otimes F &\rTo^{v_1} & \cone(\id\otimes f) &\rTo&\Sigma(\D G\otimes E)\\
\dTo^{\id}  &        & \dTo^{\D g\otimes\id}  &              &   \dDashto_{k_1}   &        &\dTo_{\id}\\
\D G\otimes E&\rTo_{\D g\otimes f} & \D F\otimes F&\rTo_{k_2} &W&\rTo&\Sigma(\D G\otimes E)
\end{diagram}

A choice for the completion morphism is

\[k_1':= \left(
\begin{array}{cc}
\id &      0          \\
0   & \D g\otimes\id
\end{array}
\right) :\cone(\id\otimes f)\lto W\]
Moreover,
\[0\to \D G\otimes E\stackrel{\id\otimes f}{\lto}\D G\otimes F\stackrel{\id\otimes g}{\lto}\D G\otimes G\to0\]
 is a semi-split exact sequence, so the morphism
\[v:=(0\qquad 1\otimes g):\cone(\id\otimes f)\lto \D G\otimes G\]
is a homotopy equivalence \ie, an isomorphism in $\BK(\CA).$ Let us denote by $\lambda$ its inverse. With  $k_1:=k_1'\lambda$, we get:
\begin{eqnarray*}
\overline{t}k_1 & = &\overline{t}k_1'\lambda\\
                            & = &\left(
\begin{array}{cc}
0 &      \epsilon
\end{array} \right)
\left(
\begin{array}{cc}
\id &      0          \\
0   & \D g\otimes\id
\end{array}\right)\lambda \\
 & = & t_G v\lambda\\
 & \simeq & t_G.
\end{eqnarray*}
Analogously, we take
\[k_3':= \left(
\begin{array}{cc}
\id &      0          \\
0   & \id\otimes f
\end{array}
\right) :\cone(\D g\otimes\id)\lto W.\]

Denoting by $\lambda'$ the inverse in $\BK(\CA)$ of
\[v'=(0\qquad\D f\otimes\id):\cone(\D g\otimes\id)\lto \D E\otimes E\]
and define $k_3:=k_3'\lambda'.$ It follows that $\overline{t}k_3\simeq t_E$ in $\BK(\CA).$
Note that $k_1=k_1'\lambda$, so $k_1 v=k_1'\lambda v\simeq k_1'$, and then
\begin{eqnarray*}
k_2(\D g\otimes\id) & =    &  k_1'v_1\\
                             &\sim & k_1 v v_1\\
      & = & k_1(\id\otimes g)
\end{eqnarray*}
that ensures the commutativity of left square on diagram \ref{diag (TC5)(a)}.
Finally, $k_2(\id\otimes f)\simeq k_3(\D f\otimes\id)$, therefore the the morphism $\overline{t}$ as previously defined makes the whole diagram commutative in $\BK(\CA).$

Moreover, $k_1$ and $k_3$ are possible choices made in Axioms (TC3) and (TC4).In fact, both come from completion to an octahedron, because, if we construct $q_1'$ and $q_1$ in the same way we did for the $k_i$'s, the  following diagram

\begin{diagram}
\cone(\D g\otimes\id) & \rTo^{v'} & \D E\otimes F & \lTo^{\D f\otimes\id}  &  DF\otimes F\\
&\luTo_{q_1'} &\uTo_{q_1} &
 \ldTo_{k_2} &\\
   \dTo^{\lambda(\D\Sigma^{-1}h\otimes g)v'}  &                    & W             &             &  \uTo_{\D g\otimes f}\\
                                                  &\ruTo^{k_1'} & \uTo_{k_1}&   \rdTo^{q_2} &    \\
\cone(\id\otimes f)  & \lTo_\lambda      & \D G\otimes G  & \rTo_{\id\otimes h} &   \D G\otimes  E
\end{diagram}
commutes.

It is clear that upper and lower triangles are commutative. Let us check that the triangle

\[\D G\otimes G\stackrel{k_1}{\lto} W \stackrel{q_1}{\lto} \D E\otimes F \stackrel{D\Sigma^{-1}h\otimes g}  {\lto}\Sigma \D G\otimes G \]
is distinguished. This is true because the following is a triangle isomorphism where the upper row is distinguished by construction:
\begin{diagram}
\cone(\id\otimes f) & \rTo^{k_1'} &W &\rTo^{q_1'} &\cone(\D g\otimes\id)&\rTo^{\gamma (\D\Sigma^{-1}h\otimes g)v'} & \Sigma\cone(\id\otimes f)  \\
 \dTo_{\lambda}&   &  \dTo_{\id} &   &\dTo^{v'} &    & \dTo_{\Sigma \lambda}\\
\D G\otimes G & \rTo_{k_1} & W &\rTo_{q_1}& \D E\otimes F & \rTo_{\D\Sigma^{-1}h\otimes g} & \Sigma \D G\otimes G
\end{diagram}
The same argument applies to $k_3$. This completes the verification of part (a) of Axiom (TC5).

Let us now treat part (b). A (TC3') type diagram for the triangles
\[\D G\stackrel{\D g}{\lto}\D F\stackrel{\D f}{\lto}DE\stackrel{\D\Sigma^{-1}h}{\lto}\Sigma \D G\]

\[E\stackrel{f}{\lto}F\stackrel{g}{\lto}G\stackrel{h}{\lto}\Sigma E\]
has its top pyramids on the form:

\begin{diagram}
\D F\otimes G           &                       & \lTo^{\id\otimes g}                  &       &\D F\otimes F        &                                & \rTo^{\D f\otimes\id}  &  &  \D E\otimes F\\
        &  \luTo^{q_3}   & \circlearrowleft  &\ldTo^{k_2}     &                & \rdTo^{k_2}  &  \circlearrowleft   &  \ruTo^{q_1} &       \\
\dTo^{\D f\otimes h}    &   +                  & W      &       +        &      \uTo^{\D g\otimes f}       &      +                                                  & W & + & \dTo_{\D\Sigma^{-1}h\otimes g} \\
                               &  \ruTo^{k_3}   &    \circlearrowleft                                                 &\rdTo^{q_2} &  &\ldTo_{q_2}       & \circlearrowleft &\luTo_{k_1} &\\
\D E\otimes E           &                       & \rTo_{\D\Sigma^{-1}h\otimes \id} & & \D G\otimes E  &         &\lTo_{\id\otimes h}& &\D G\otimes G
\end{diagram}

The dual of a (TC3) type diagram for the triangles
\[ E\stackrel{f}{\lto}F\stackrel{g}{\lto}G\stackrel{h}{\lto}\Sigma E\]

\[\D G\stackrel{\D g}{\lto}\D F\stackrel{\D f}{\lto}\D E\stackrel{\D\Sigma^{-1}h}{\lto}\Sigma \D G\]
has its top pyramids on the form:
\begin{diagram}
\D(F\otimes \D G)           &                       & \lTo^{\D(\id\otimes \D g)}                  &       &\D(F\otimes \D F)        &                        & \rTo^{\D(f\otimes\id)}  &  &  \D(E\otimes \D F)\\
 &  \luTo^{\D p_1}   &    \circlearrowleft  &\ldTo^{\D j_2}     &                & \rdTo^{\D j_2}  &  \circlearrowleft      &  \ruTo^{\D p_3} &       \\
 \dTo^{\D(f\otimes \D\Sigma^{-1}h)}  &                       & \D V &         &  \uTo^{\D(g\otimes \D f)}                                          &                     & DV &    &   \dTo_{\D(h\otimes \D g)}  \\
                               &  \ruTo^{\D j_1}   &  \circlearrowleft                                                   &\rdTo^{\D p_2} &  &\ldTo_{\D p_2}& \circlearrowleft    &\luTo_{\D j_3} &\\
\D(E\otimes \D E)           &                       & \rTo_\ & &\D(G\otimes \D E)  &         &\lTo_{\D(\id\otimes \D h)}& &\D(G\otimes \D G)
\end{diagram}

In this last diagram, morphisms $\D p_i$ and $\D j_i$ are determined by the morphisms $p_i$ and $j_i$, but as all objects considered are strongly dualizable, the reciprocal is also true. In fact, denote by $P_i$ and $J_i$ the morphisms chosen after dualizing  all objects and all morphism different from $p_i$ and $j_i$ on the (TC3) type diagram, and by $A,$ $B$ $C$ or $D$ any object appearing in any vertex. All  these objects are strongly dualizable, so, using the isomorphism $\beta: A\lto \D\D A$ from Lemma \ref{ddual}, and the fact that due to the naturality of $\beta,$  for any morphisms $u$, $v$, it holds that $\beta(\D u\otimes v)=\D\D(\D u\otimes v)\rho$. Then if we dualize again, it arises a (TC3) type diagram isomorphic to one of the same type with $V$ on the apices and $\beta^{-1}\D P_i \beta$, $\beta^{-1}\D J_i \beta$ on the edges.The triangle
\[
A\otimes B \stackrel{\beta^{-1}\D P_i\beta}{\lto} V \stackrel{\beta^{-1}\D J_i\beta}{\lto} C\otimes D \stackrel{+}{\to}
 \]
is distinguished because its double dual is.
Moreover, if $P_i$ and $J_i$ were defined from the completion to an octahedron, the same is true for $DP_i$ and $DJ_i$ and hence for $\beta^{-1}\D P_i\beta$ and $\beta^{-1}\D J_i\beta$.
Then, to achieve our result, the connection morphisms of the last diagram, $\D P_i$ and $\D J_i$,  can be chosen after dualizing.

Returning to the construction, it must be found
 an isomorphism between the chained octahedron  of the (TC3') type diagram for the triangles with successive morphisms $(\D g,\D f,\D\Sigma^{-1}h)$ and $(f,g,h)$ and the chained octahedron of the (TC3) type diagram for $(f,g,h)$ and $(\D g,\D f,\D\Sigma^{-1}h)$.

We chose $W:=\cone(\D g\otimes f)$ with $k_2$ and $q_2$ the canonical projection and injection morphisms. Looking at the diagram,

\begin{diagram}
\D G\otimes E & \rTo^{\D g\otimes f} &\D F\otimes F &\rTo^{k_2} & W &\rTo^{q_2} & \Sigma(\D G\otimes E)  \\
 \dTo_{\xi}&   &  \dTo_{\xi} &   & \dDashto_{\overline{\xi}}    &    & \dTo_{\Sigma \xi}\\
\D(G\otimes \D E) & \rTo_{\D(g\otimes \D f)} & \D(F\otimes \D F) &\rTo_{J_2}& \D V & \rTo_{P_2} & \Sigma (\D(G\otimes \D G))
\end{diagram}

Lemma \ref{funxi} ensures that the firs square is commutative. The lower triangle is distinguished, because it is the dual of a distinguished one, so we can complete to a morphism of triangles. A possible choice of the remaining morphism is
\[\overline{\xi}:=\left(\begin{array}{cc} \xi & 0\\
0 &\xi      \end{array}\right)\]
see \cite{Lip} 1.4.3. Now, to define the morphism of octahedra we were looking for, we chose for each vertex different from  $W$  the corresponding map $\xi$ and $\overline{\xi}:W\rightarrow \D V$.
 Using Lemma \ref{funxi}, all diagrams arising where it does not occur any $P_i$, $J_i$, $k_i$ or $q_i$ commute. Then, to prove that it is actually  a morphism of octahedra, we are reduced to check the commutativity when some of these arrows occur. To do so, let's consider the diagram
\begin{diagram}
  &                                     & \D(G\otimes \D E)&                                  & \\
 & \ldTo^{\D(\id\otimes \D f)}&                          &\rdTo^{\D(g\otimes\id)}& \\
\D(G\otimes \D F)&\rTo_{\D(g\otimes\id)}&\D(F\otimes \D F)& \lTo_{\D(\id\otimes \D f)} &\D(F\otimes \D E)\\
\dTo^{\D(\id\otimes \D g)} & & \dTo_{J_2} & & \dTo_{\D(f\otimes\id)}\\
\D(G\otimes \D G)&\rTo_{J_3}&\D V& \lTo_{J_1} &\D(E\otimes \D E)\\
\dTo^{\phi}                   &\ruTo_{J_3'} & &\luTo_{J_1'}&\dTo_{\phi'}\\
\cone(\D(\id\otimes \D f))&                   &  &                  &\cone(\D(g\otimes\id))
\end{diagram}
where $\phi$ and $\phi'$ are homotopy equivalences. The definition of the $J_i$'s is as follows:  completing the diagram
\begin{diagram}
\D(G\otimes \D E) & \rTo^{\D(g\otimes\id)} &\D(F\otimes \D E) &\rTo^{(1,0)} & \cone(\D(g\otimes\id)) &\rTo^{\text{{\scriptsize$\left(\begin{array}{c} 1\\ 0  \end{array}\right)$}}} & \Sigma(\D G\otimes E)  \\
 \dTo_{\id}&   &  \dTo_{\D(\id\otimes \D f)} &   & \dDashto_{J'_1}    &    & \dTo_{\id}\\
\D(G\otimes \D E) & \rTo_{\D(g\otimes \D f)} & \D(F\otimes \D F) &\rTo_{J_2}& \D V & \rTo_{P_2} & \Sigma (\D(G\otimes \D G))
\end{diagram}
to a morphism of standard tringles. We may take \[J'_1:=\text{\scriptsize$\left(\begin{array}{cc} 1 & 0\\
0 & \D(\id\otimes \D f) \end{array}\right)$}\]

In an analogous way we define $J'_3=\cone(\D(id\otimes \D f))\rightarrow \D V$ as
\[J'_3:=\text{\scriptsize$\left(\begin{array}{cc} 1 & 0\\
0 & \D(\id\otimes \D f) \end{array}\right)$}\]

Finally, define:
$J_1:=J_1' \phi'$ and $ J_3=J'_3\phi$ completing the definition of the remaining morphisms in the diagram We proceed now to check the  commutativity of the squares  required to have a morphism between chained octahedra.

We have to check that $\overline{\xi}k_3=\xi J_1.$ To get it, we will define maps $\lambda'$ and $\phi'$ such that $k_3=k'_3\lambda'$ and $J_1=J'_1\phi',$ so we have two squares

\begin{diagram}
\D E\otimes E & \rTo^{\lambda'} &\cone(\D g\otimes\id) &\rTo^{k_3'} & W \\
 \dTo_{\xi}&   &  \dTo_{\text{\scriptsize $\left(\begin{array}{cc} \xi & 0 \\ 0 & \xi \end{array}    \right)$}} &   &\dTo_{\text{ \scriptsize $\left(\begin{array}{cc} \xi & 0 \\ 0 & \xi \end{array}    \right)$}}\\
\D(E\otimes \D E) & \rTo_{\phi'} & \cone(\D(g\otimes\id)) &\rTo_{J_1'}& \D V
\end{diagram}

Note that the right square commutes, because $\xi(\id\otimes f)=\D(\id\otimes Df)\xi$ by Lemma \ref{funxi}. Let us build the left one. We assume that there is a semi-split short exact sequence
\[0\lto E\stackrel{f}{\lto}F\stackrel{g}{\lto}G\lto0\]
associated to the triangle. Therefore,
\[0\lto \D G\stackrel{\D g}{\lto}\D F\stackrel{\D f}{\lto}\D E\lto0\]
is also a semi-split short exact sequence.
We choose splittings  $\psi$, $\varphi$ for $f$ and $g$ respectively, so  $\D\psi$ and $\D\varphi$ are splittings for $\D f$ and $\D g$ respectively.

Looking at \cite[Example (1.4.3)]{Lip}  we put
\[\lambda':=\text{\scriptsize$\left(\begin{array}{c} \D\psi\otimes id\\
 a \end{array}\right)$}\text{ and } \phi':=\text{\scriptsize$\left(\begin{array}{c} \D(\psi\otimes id)\\
 b \end{array}\right)$}\]
where the maps of complexes $a$ and $b$ are defined by

\[a^n=(D(\psi^{n+1})\otimes\id)d^n_{\D E\otimes E} -d^n_{\D F\otimes E}(D(\psi^n)\otimes\id))\]
\[b^n=\D(\psi^{n+1}\otimes\id)d^n_{\D(E\otimes\D E)}  - d^n_{\D(F\otimes \D E)}(\D(\psi^n\otimes\id))\]
for $n\in\mathbb{Z}$. A computation using Lemma \ref{funxi} shows that

\[\text{\scriptsize  $\left(\begin{array}{cc} \xi & 0
\\ 0 & \xi \end{array}\right)$}\lambda'=\text{\scriptsize $ \left(\begin{array}{cc} \xi & 0 \\ 0 & \xi \end{array}\right)\left(\begin{array}{c}\D \psi\otimes id\\
a\end{array}\right)$}=\text{\scriptsize $\left(\begin{array}{c}\D(\psi\otimes id)\\
b\end{array}\right)$}\xi=\phi'\xi\]

Following a similar argument we see that $overline{\xi}k_1=\xi J_3.$

There is left to prove the commutativity of the following squares:

\begin{diagram}\label{diagrams}
W &\rTo^{q_i} &\D F\otimes G\\
\dTo^{\overline{\xi}} & & \dTo_{\xi}\\
\D V &\rTo_{P_j} & \D(F\otimes \D G)
\end{diagram}
for $(i,j)=(1,3)$ or $(3,1)$. We argue as before. To complete the diagram

\begin{diagram}
\D G\otimes E & \rTo^{\D g\otimes f} &\D F\otimes F &\rTo^{k_2} & W&\rTo^{q_2}&\Sigma(\D G\otimes E)\\
\dTo_{\D g\otimes\id}&   &\dTo_{id} &   &\dDashto_{q'_i}&   &    \dTo_{\Sigma(\D g\otimes\id)}\\
\D F\otimes E & \rTo_{\id\otimes f)} &\D F\otimes F &\rTo_{} &\cone(\id\otimes f)&\rTo&\Sigma(DF\otimes E)
\end{diagram}
to a morphism of standard triangles. We may take

\[q_3'=\text{\scriptsize $\left(\begin{array}{cc} \D g\otimes\id & 0 \\ 0 & \id \end{array}\right)$}\]
 and
$q_i:=\chi q_i'$
where $\chi:\cone(\id\otimes f)\to \D F\otimes G$ is the homotopy equivalence defined as  $\chi:=(\id\otimes f\quad 0)$.
Similarly, we take $P'_j:\D V\rightarrow\cone(\D(id\otimes\D f))$ as

\[P'_j=\text{\scriptsize$\left(\begin{array}{cc} \D(g\otimes\id) & 0 \\ 0 & \id \end{array}\right)$}\]
 and define
$P_j:=\chi' P_j'$
where $\chi':=\left(\D(id\otimes\D g)\quad 0 \right)$ is the obvious homotopy equivalence. To prove the desired commutativity, decompose the diagram \ref{diagramas} as
\begin{diagram}
W & \rTo^{q_3'} &\cone(\id\otimes f) &\rTo^{\chi} & \D F\otimes G \\
 \dTo^{\overline{\xi}}&      &\dTo^{\text{{\scriptsize$\left(\begin{array}{cc} \xi & 0 \\ 0 & \xi \end{array}    \right)$}}}& &\dTo_{\xi}\\
\D V & \rTo_{P_1'} & \cone(D(\id\otimes \D f)) &\rTo_{\chi'}& \D(F\otimes \D G).
\end{diagram}
Both squares commute by a further use of Lemma \ref{funxi}.

Then, we have defined a morphism between three chained octahedra, because all diagrams arising are commutative. Note also that the $P_i$'s and $J_i$'s morphisms come from completion to an octahedron, because they are constructed in the same way as $q_i$ and $k_i$, which come from such a completion as sawn in (TC5)(a).

Finally, recall that  an involution of the later (TC3) type diagram is just a (TC3) type diagram for the triangles $(Dg,Df,D\Sigma^{-1}h)$, $(f,g,h)$ whose connection morphism come from completion to an octahedron, so the condition (TC4)  is satisfied for a (TC3') type diagram and the mentioned involution, just because we are in the hypothesis of the Axiom. \qed.

As is pointed out May, there is a dual version of this Axiom. We will state it for latter use.

\begin{corollary}\label{TC5dual}
 Given a (TC3) type diagram for triangles $(f,g, h)$ and $(\D g,\D f,\D\Sigma^{-1}h)$ as considered in (TC5)(b), there is a map $\overline{u}:S\to V$ making  the diagram
\begin{diagram}
 &                                          &    S                 &                        &  \\
 & \ldTo^{(u_G,u_E)} &\dTo\overline{u}&\rdTo^{u_F}& \\
(\D G\otimes G)\oplus(\D E\otimes E) & \lTo_{(j_1,j_3)} & V &\rTo_{j_2}& \D F\otimes F
\end{diagram}
commutative, where by $u$ we denote the map induced by the unit of the tensor $-$ hom adjuntion.
\end{corollary}

\proof It follows from \cite[Lemma 4.14 (TC5)(a')]{May}. \qed

\begin{remark}
Note that the verification of any of the compatibility Axioms separately  is trivial; we can define inmediatenly morphisms with the required properties. The Axioms make sense when we consider them all together; the problem is that in each Axiom we must make a choice, for an object and for certain morphisms, but Axioms $(TC4)$ and both parts of $(TC5)$ refer to the choice we made before  when using the third Axiom of triangulated categories to get a non canonical arrow between the third objects on distinguished triangles. What we do here is to give explicitely such a good choice, that is in some sense the expected one.
\end{remark}

\section{Additivity of the trace of an endomorphism}

One of the features of abstract duality as defined in \cite[\S III.1]{LMS} is that an endomorphism of a strongly dualizable object possesses a trace, that is a well-defined element of the ring of endomorphisms of the unit object. We will adopt a more general view point allowing for an orientation or genus that will give us the correct formalism for studying cohomological characters.

To fix ideas, let $\mathbf{T}$ be a triangulated category with a compatible closed structure. An orientation with values in an object $C\in\mathbf{T}$ is a natural transformation
\[A:id_\mathbf{T}\rightarrow -\otimes C.\]
It provides an {\it orientation map} for every object $E\in\mathbf{T}$
\[A_E:E\rightarrow E\otimes C.\]

Recall that, by definition $[E,-]\cong \mathrm{D}E\otimes-$ therefore the unit and counit of the tensor-hom adjunction induce the following natural maps
\[u_E:S\rightarrow E\otimes\mathrm{D}E\text{ and }t_E:\mathrm{D}E\otimes E\rightarrow S.\]

Given an endomorphism of a strongly dualizable object $\phi:E\rightarrow E$ we are going to define its Lefschetz invariant with respect to the orientation $A$ as the following composition
\[S\stackrel{u}{\to}E\otimes\mathrm{D}E\stackrel{\gamma}{\to}\mathrm{D}E\otimes E\stackrel{id\otimes\phi}{\lto}\mathrm{D}E\otimes E\stackrel{id\otimes A_E}{\lto}\mathrm{D}E\otimes E\otimes C\stackrel{t\otimes id}{\lto}S\otimes C\cong C\]
we will denote it by $Lef(\phi,A)\in\Hom_\mathbf{T}(S,C).$ Note that in \cite{May} itisused the notation $\tau_A(\phi).$

If the orientation is trivial \ie $C=S$ and $A=\rho^{-1}$ given by the natural isomophism $E\cong E\otimes S$ then we denote it by $Tr(\phi):=Lef(\phi,\rho^{-1}).$ For the rest of the section we will deal with the case of a trivial orientation and $\mathbf{T}=\mathbb{D}(\mathcal{A})$ where $\CA$ is a category satisfying the hypothesis of the first section.

We start with the following useful homological algebra result. In short, it says that we an replace the objects and maps in a homotopy-commutative morphism of triangles by an isomorphic one in the homotopy category such that it is commutative already in the category of complexes.

\begin{lemma}\label{lemacyl} Let $\mathbf{\mathcal{C}}$ be an additive category and $\mathbb{K}(\mathcal{C})$ its homotopy category. Given a distinguished triangle $E\stackrel{f}{\lto}F\stackrel{G}{\lto}G\stackrel{h}{\lto}\Sigma E$ and a homotopy commutative  diagram
\begin{diagram}
E               & \rTo^{f} & F  & \rTo^{g} & G  & \rTo^{h} & \Sigma E \\
\dTo^{\phi} &             & \dTo_{\psi}   &      &   \dTo_{\omega}       &     & \dTo_{\Sigma\phi}\\
E               & \rTo_{f} & F  &\rTo_{g} &  G & \rTo_{h}  & \Sigma E
\end{diagram}
it is possible to replace the objects and the maps in the diagram in such a way that we obtain a morphism of triangles isomorphic to the previous one (in $\mathbb{K}(\mathcal{C})$) but whose underlying diagram in $\mathbf{C}(\mathcal{C})$ is commutative.
 \end{lemma}

\proof
Note that all the squares in the diagram are commutative in $\mathbb{K}(\mathcal{C}),$ \ie are commutative up to homotopy. A first remedy is to replace $g$ by $\cone(f),$ obtaining a new diagram
\begin{diagram}
E               & \rTo^{f} & F  & \rTo^{v} & \cone(f)  & \rTo^{p} & \Sigma E \\
\dTo^{\phi} &    \textcircled{\scriptsize{1}}         & \dTo_{\psi} & \textcircled{\scriptsize{2}} &\dTo_{\omega'}   &   \textcircled{\scriptsize{3}}    &      \dTo_{\Sigma\phi}\\
E               & \rTo_{f} & F  &\rTo_{v} &  \cone(f) & \rTo_{p}  & \Sigma E
\end{diagram}
where we may choose
\[\omega':=\left(\begin{array}{cc} \Sigma\phi & 0 \\
                                     s        & \psi  \end{array} \right)\]
with $s:E\rightarrow\Sigma^{-1}F$ the homotopy such that $f\phi-\psi f= s d_E+d_F s.$ Note that $\omega'$ is a morphism o complexes and that the second and third  squares are now commutative in $\mathbf{C}(\mathcal{C}).$ Also, it is clear that $\cone(f)$ is isomorphic to $G$ in $\mathbb{K}(\mathcal{C}).$ To make commutative the first square we must do a further replacement. Take $\cyl(f)=\Sigma^{-1}\cone(h')$ with differential
\[d_{\cyl(f)}=\left(\begin{array}{ccc}
 -d_E & 0 & 0\\
f & d_F  & 0 \\
-\id & 0 & -d_E
\end{array}\right).\]
This gives us a new diagram
\begin{diagram}
\cone(f)               & \rTo^{h'} & \Sigma E  & \rTo^{f'} & \cyl(f)  & \rTo^{g''} & \Sigma\cone(f) \\
\dTo^{\omega'} &  \textcircled{\scriptsize{3}}  & \dTo_{\Sigma\phi} & \textcircled{\scriptsize{${1}^\prime$}} &\dTo_{\psi'}   & \textcircled{\scriptsize{${2}^\prime$}} &      \dTo_{\Sigma\omega'}\\
\cone(f)               & \rTo_{h'} & \Sigma E  &\rTo_{f'} &  \cyl(f) & \rTo_{g''}  & \Sigma\cone(f)
\end{diagram}
where we may choose
\[\psi':=\left(\begin{array}{ccc} \Sigma\phi & 0 &0\\
s & \psi & 0 \\
0 & 0 &\phi
\end{array}\right)\]

Note that $\psi'$ is a morphism of complexes and that the second and third squares in the preceding diagram are commutative in $\mathbf{C}(\mathcal{C}).$ Also, by its definition, $\cyl(f)$ is isomorphic to $\Sigma F$ in $\mathbb{K}(\mathcal{C}).$ Desuspending the above morphism of triangles we reach the desired conclusion.\qed

\begin{remark}
Intuitively speaking, the lemma says that we can codify the information given by a homotopy between morphisms deforming the domain and target objects, in such a way that the information is now contained on these new objects, more precisely, on the differentials of complexes. This is specially clear on topological examples.
\end{remark}

\begin{theorem}\label{addtraces}{\bf Additivity of traces}.
Let $E\stackrel{f}{\to}F\stackrel{g}{\to}G\stackrel{h}{\to}\Sigma E$ be a distinguished triangle of strongly dualizable objects and consider the following commutative diagram of solid arrows
\begin{diagram}
E               & \rTo^{f} & F  & \rTo^{g} & G  & \rTo^{h} & \Sigma E \\
\dTo^{\phi} &             & \dTo_{\psi}   &      & \dDashto_{\omega}         &     & \dTo_{\Sigma\phi}\\
E               & \rTo_{f} & F  &\rTo_{g} &  G & \rTo_{h}  & \Sigma E
\end{diagram}
There is at least an explicit choice of a morphism  $\omega:G\to G$ that makes the whole diagram a morphism of triangles, and such that
\[Tr(\psi)=Tr(\phi)+Tr(\omega).\]
\end{theorem}
\proof
By taking $q$-injective resolutions in $\BD(\CA)$ we are at once reduced to the case where the left square is commutative in $\BK(\CA)$ with commuting homotopy $s:E\rightarrow\Sigma^{-1}F.$ Using the lemma, and denoting $G':=\cone(f)$ $F'=\cyl(f),$ it suffices to deal with distinguished triangles

\[E\stackrel{f'}{\lto}F'\stackrel{g''}{\lto}G'\stackrel{h'}{\lto}\Sigma E.\]
and
\[\mathrm{D}G'\stackrel{-\Sigma^{-1}\mathrm{D}h'}{\lto}\mathrm{D}F'\stackrel{\mathrm{D}g''}{\lto}\mathrm{D}E\stackrel{\mathrm{D}f'}{\lto}\Sigma \mathrm{D}G'\]
and replace $\psi$ and $\omega$ by $\psi'$ and $\omega'$ respectively.

The additivity formula follows from the commutativity of the outer maps of the following diagram

\begin{diagram}
                                &                          &  S                     &                    &    \\
                                &\ldTo^{u_{G'},u_E}         & \dTo_{\overline{u}} &\rdTo^{u_{F'}}  &      \\
  (G'\otimes \mathrm{D}G')\oplus(E\otimes \mathrm{D}E) &\lTo_{(j_3,j_1)}  &V &\rTo_{j_2}    & F'\otimes DF'\\
\dTo^{(\gamma,\gamma)}&      \mathbf{\mathrm{nat}}                       & \dTo_{\gamma} &  \mathbf{\mathrm{nat}}   &\dTo_{\gamma} \\
                                 &                          & \overline{V}    &                    &    \\
                                     &\ldTo^{(\overline{j_1},\overline{j_3})} &  &\rdTo^{\overline{j_2}} & \\
(\mathrm{D}G'\otimes G')\oplus(\mathrm{D}E\otimes E) &  &\star &    &\mathrm{D}F'\otimes F' \\
\dTo^{(\id\otimes\omega',\id\otimes\phi)}&\rdTo_{(\overline{k_1},\overline{k_3})}                              &   &\ldTo_{\overline{k_2}} &    \dTo_{\id\otimes\psi'} \\
                                       & \textcircled{\scriptsize{1}} & \overline{W}&\textcircled{\scriptsize{2}} & \\
         &           &  \dDashto_{m} &    &    &  \\
(\mathrm{D}G'\otimes G')\oplus(\mathrm{D}E\otimes E) &\rTo^{(\overline{k_1},\overline{k_3})}  & \overline{W}&\lTo^{k_2}    &\mathrm{D}F' \otimes F'\\
                                 &\rdTo_{(t_{G'},t_E)} & \dTo_{\overline{t}}              &\ldTo_{t_{F'}} &    \\
                                &                          &  S &                    &
\end{diagram}

To check this commutativity, we decompose the diagram into smaller subdiagrams. Let us see why they do commute. The upper triangles commute by Corollary \ref{TC5dual}, while the lower ones by Axiom (TC5)(a), as proved in Theorem \ref{tc5h}. For the middle part, the  trapezoids marked by $\mathbf{\mathrm{nat}}$ commute by \ref{barra} and the central rhomboid marked by $\star$ commutes by Axiom (TC5)(b) after applying an involution (via $\gamma$) which is possible in view of \ref{barra} again.

Therefore, to finish the proof we need to find an endomorphism $m$ of $\overline{W}$ such that subdiagrams \text{\textcircled{\scriptsize{1}}} and \text{\textcircled{\scriptsize{2}}} commute. Consider the following diagram of distinguished triangles

\begin{diagram}
(\mathrm{D}G'\otimes E) & \rTo^{\mathrm{D}g''\otimes f'} & \mathrm{D}F'\otimes F' & \rTo^{k_2} & \overline{W} &\rTo & \Sigma(\mathrm{D}G'\otimes E) \\
\dTo^{(\id\otimes\phi)} & & \dTo_{\id\otimes\psi'} & & \dDashto_m    & &\dTo_{\Sigma (\id\otimes\phi)}\\
(\mathrm{D}G'\otimes E) & \rTo_{\mathrm{D}g''\otimes f'} & \mathrm{D}F'\otimes F' & \lTo_{k_2} & \overline{W} & \rTo & \Sigma(\mathrm{D}G'\otimes E)
\end{diagram}
and, as the first square commutes in $\mathbf{C}(\CA)$ we can choose:

\[m:=\left(\begin{array}{cc} \id_{\mathrm{D}G'}\otimes\phi & 0 \\
                                           0                              & \id_{\mathrm{D}F'}\otimes\psi'
 \end{array}\right)\]
making our second square commutative in $\mathbf{C}(\CA)$ (and in $\BK(\CA)$), but this is exactly the lower right trapezoid. Now to
 prove the commutativity of the remaining one it is enough to check that the diagram
\begin{diagram}\label{circ2}
(\mathrm{D}G'\otimes G')                & \rTo^{k_3} & \overline{W}   & \lTo^{k_1} & \mathrm{D}E\otimes E  \\
\dTo^{\id\otimes\omega'} &                 & \dTo_{m} & & \dTo_{\id\otimes\phi}\\
(\mathrm{D}G'\otimes G')                & \rTo_{k_3} & \overline{W}    & \lTo_{k_1} & \mathrm{D}E\otimes E.
\end{diagram}

commutes. Now that $m$ is fixed (even as a map of complexes) we may make a further change of distinguished triangles within the isomorphism class in $\BK(\CA).$ Replace $(f',g'',h')$ by
\[E\stackrel{f'}{\lto}F'\stackrel{g'''}{\lto}G''\stackrel{h''}{\lto}\Sigma E\]
where $G'':=\cone(f')$ and $g'''$ and $h''$ are the canonical morphisms.
We need to define an endomorphism $\omega'':G''\to G''$ compatible with $\omega',$ namely
\[\omega''=\left(\begin{array}{cc}\Sigma\phi & 0 \\
                            0       \ \psi' \end{array}\right).\]
Let us make this explicit. We are going to define an isomorphism $\theta:G'\to G''$ in $\BK(\CA)$ such that the following square
\begin{diagram}
G'&\rTo^{\theta}&G''\\
\dTo^{\omega'}& & \dTo_{\omega''}\\
G'&\rTo^{\theta}&G''
\end{diagram}
commutes. To this end express $G'\stackrel{gr}{=}\Sigma E\oplus F$ and $G''\stackrel{gr}{=}\Sigma E\oplus F',$ where "$\stackrel{gr}{=}$" means equality as graded objects. Define morphisms $\alpha_1:F\to F'$ and $\alpha_2:\Sigma E\to F'$ by
\[\alpha_1:=\left(\begin{array}{ccc} \id & 0 & 0\end{array}\right)\text{ and }\alpha_2:=\left(\begin{array}{ccc} 0 & \id & 0\end{array}\right)\]
respectively, representing $F'\stackrel{gr}{=}\Sigma E\oplus F\oplus E.$ Note that $\mathbf{\alpha}:=\left(\begin{array}{cc} \alpha_1& \alpha_2\end{array}\right)$ is a section of $g''.$ The morphism
\[\theta:=\left(\begin{array}{cc} \id        & 0 \\
                                    \alpha_1 & \alpha_2 \end{array}\right)\]
is a homotopy isomorphism as it follows the pattern of \cite[Example 1.4.3]{Lip}. With these definitions,
\[\theta\omega'=\left(\begin{array}{cc} \id        & 0 \\
                                    \alpha_1 & \alpha_2 \end{array}\right)\left(\begin{array}{cc} \Sigma\phi        & 0 \\
                                    s & \psi \end{array}\right)=
                                     \left(\begin{array}{cc} \Sigma\phi        & 0 \\
                                    \alpha_1\Sigma\phi+\alpha_2 s & \alpha_2\psi' \end{array}\right)\]
and
\[\omega''\theta=\left(\begin{array}{cc} \Sigma\phi        & 0 \\
                                    s & \psi' \end{array}\right)
                                    \left(\begin{array}{cc} \id        & 0 \\
                                    \alpha_1 & \alpha_2 \end{array}\right)=
                                    \left(\begin{array}{cc} \Sigma\phi        & 0 \\
                                    \alpha_1\psi' & \alpha_2\psi\end{array} \right).\]
Note first that
\[\alpha_1\Sigma\phi+\alpha_2 s=\left(\begin{array}{ccc} \Sigma\phi& s & 0 \end{array}\right)\]
and
\[\alpha_1\psi'=\left(\begin{array}{ccc}\id & 0 & 0\end{array}\right)
\left(\begin{array}{ccc}
\Sigma\phi & 0   & 0\\
s          &\psi & 0\\
0          &  0  & \phi
\end{array}\right)=
\left(\begin{array}{ccc} \Sigma\phi & s& 0\end{array}\right).\]
For the remaining entry, note that
\[\alpha_2\psi'=\left(\begin{array}{ccc}0 & \id & 0\end{array}\right)
\left(\begin{array}{ccc}
\Sigma\phi & 0   & 0\\
s          &\psi & 0\\
0          &  0  & \phi
\end{array}\right)=
\left(\begin{array}{ccc} 0 & \psi & 0\end{array}\right),\]
\[\alpha_2\psi=\left(\begin{array}{ccc} 0 & \psi & 0\end{array}\right).\]

Let us prove now the commutativity of the left part of diagram \ref{circ2}. Consider the following decomposition

\begin{diagram}
\mathrm{D}G'\otimes G'  & \rTo^{\id\otimes\theta} & \mathrm{D}G'\otimes G''   & \rTo^{k_3'} & \overline{W}  \\
\dTo^{\id\otimes\phi} &                 & \dTo^{\id\otimes\omega''} & &\dTo_m\\ \mathrm{D}G'\otimes G'  & \rTo_{\id\otimes\theta} & \mathrm{D}G'\otimes G''   & \rTo^{k_3'} & \overline{W}
\end{diagram}
The corresponding left part commutes by the defining property of $\theta$ and the right part because we may take
\[k'_3=\left(\begin{array}{cc}
\id_{\mathrm{D}G'}\ot\id_{\Sigma E} & 0 \\
 0                                 & \mathrm{D}g'\otimes\id_{F'}
\end{array}\right)\]
so both paths yield
\[mk'_3=\left(\begin{array}{cc}
\id_{\mathrm{D}G'}\otimes\Sigma\phi  & 0 \\
 0                                   & \mathrm{D}g'\otimes\psi'
\end{array}\right)=k'_3(\id\otimes\omega'').\]

To finish the proof, let us prove the commutativity of the right part of diagram \ref{circ2}. Describe, as in the proof of Theorem \ref{tc5h} $\overline{W}=\cone(\mathrm{D}(g')\otimes f').$ Recall that the cone map $\mathrm{D}(g'):\mathrm{D}(G')\to\mathrm{D}F$ is homotopically equivalent t $\mathrm{D}E.$ Denote by $\Phi:\mathrm{D}E\to\cone(\mathrm{D}(g'))$ an isomorphism in $\BK(\CA)$ constructed as in \cite[Example 1.4.3]{Lip}. We may decompose the square as
\begin{diagram}
\mathrm{D}E\otimes E & \rTo^{\Phi\otimes\id} & \cone(\mathrm{D}g')\otimes E   & \rTo^{k_1'} & \overline{W}  \\
\dTo^{\id\otimes\phi} &                 & \dTo^{\id\otimes\phi} & &\dTo_m\\ \mathrm{D}E\otimes E & \rTo^{\Phi\otimes\id} & \cone(\mathrm{D}g')\otimes E   & \rTo^{k_1'} & \overline{W}
\end{diagram}
Note that the left subdiagram is obviously commutative. For the remaining part, $k_1$ may be factored as $k'_1(\Phi\otimes\id)$ with
\[k'_1=\left(\begin{array}{cc}
\id  & 0 \\
 0   & \id\otimes f'
\end{array}\right)\]
therefore the left subdiagram commutes just because $\psi' f'=f'\phi.$
\qed

\section{Additivty of the Chern Character}

From now on, we will stick with  part (iii) of Example \ref{exemplos} in the first section. Consider a quasi-compact and separated scheme $X.$  We will assume that $X$ is smooth over a base scheme $S$ which is connected. In this case the dimension of fibers is constant and we will denote it by $n$. In this case, our category $\CA:=\CA_{qc}(X)$ will be the derived category of quasi-coherent sheaves on $X$. As we have already remarked, its derived category $\mathbb{D}(\CA_{qc}(X))$ satisfies all our hypothesis, therefore, by Theorem \ref{th_ct} is a triangulated category with a compatible closed structure. The closed structure is the usual one and the unit object is $S=\CO_X$.

We are going to define an orientation in $\mathbb{D}(\CA_{qc}(X))$ in the sense of the previous section, through the Atiyah class. This was introduced by O'Brian, Toledo and Tong in \cite{OTT}. We will give here a treatment adapted to our context, following \cite[\S 16]{EGA IV}.

\begin{paragrafo}\bf{Principal parts.}
\end{paragrafo}
Consider the diagonal embedding
\[\delta:X\lto X\times_S X.\]

Denote by $X_\delta$ the image scheme and by $X_\delta^{(1)}$ its first infinitesimal neighborhood, \ie the only subscheme of $X\times_S X$ defined by the sheaf of ideals $\CI^2$ where $\CI$ is the sheaf of ideals of $X_\delta$. Denote, for $i\in\{i,2\}$ by $p_i:X\otimes_SX\to X$ the canonical projections and by $h:X_\delta^{(1)}\to X\times_SX$ the canonical embedding. Let $p_1^{(1)}:=p_i h$.

Let $\CF\in\CA_{qc}(X)$ and define the sheaf of 1-principal parts of $\CF$ as
\[\CP^1_{X|S}(\CF):=(p_1^{(1)})_*(p_2^{(1)})^*\CF\]
We will abbreviate $\CP^1_{X|S}(\CO_X) $ by $\CP^1_{X|S}$. There is an obvious counit map
\[\CP_{X|S}^1(\CF)\lto\CF\]
Let $\Omega_{X|S}^1$ or simply $\Omega_X^1$ denote the sheaf of differentials. A local computation yields an exact sequence
\begin{equation}\label{ppseq}
0\lto\Omega_{X|S}^1\lto\CP_{X|S}^1\lto\CO_X\lto0
\end{equation}
split as sequence of sheaves of abelian groups, but not a s a sequence of $\CO_X$-modules (the splitting, is in fact $x^{-1}\CO_S$-linear, where $x:X\to S$ denotes the structure map. Due to this splitting and the isomorphism $\CP_{X|S}^1(\CF)\cong\CP_{X|S}^1\otimes\CF$ from \cite[(16.7.2.1)]{EGA IV}), we have for every $\CF$ an exact sequence
\[0\lto\Omega_{X|S}^1\otimes\CF\lto\CP_{X|S}^1(\CF)\lto\CF\lto0\]
obviously functorial in $\CF$.

\begin{paragrafo}\bf{The Atiyah class.}
\end{paragrafo}
Recall that an object $\CE\in\mathbb{D}(\CA_{qc}(X))$ is strongly dualizable if and only if it is perfect \cite{Neeperf}. We say that $\CE$ is a perfect complex if for every $P\in X$ there is an open neighborhood $P\in U\subset X$ such that, denoting $j:U\to X$ the canonical inclusion, the complex $j^*(\CE)$ is quasi-isomorphic to a bounded complex made up of locally free finite-type Modules over $U$.

Take a specific representative $\CE$ in $\mathbf{C}(\CA_{qc}(X))$. Tensor it with the sequence \ref{ppseq}. We obtain an exact sequence of complexes
\begin{equation}\label{ppseqder}
0\lto\CE\otimes\Omega_{X|S}^1\lto\CE\otimes\CP_{X|S}^1\lto\CE\lto0
\end{equation}
It corresponds in $\mathbb{D}(\CA_{qc}(X))$ to a triangle
\[\CE\otimes\Omega_{X|S}^1\lto\CE\otimes\CP_{X|S}^1\lto\CE\stackrel{\at_\CE^1}
{\lto}\CE\otimes\Sigma\Omega_{X|S}^1\]
In other words,
\[\at_\CE^1\in\Hom_{\BD(\CA)}(\CE,\CE\otimes\Sigma\Omega_{X|S}^1)=
Ext^1(\CE,\CE\otimes\Sigma\Omega_{X|S}^1)\]
is the class of the extension (\ref{ppseqder}).

\begin{remark} The Atiyah class is an invariant that depends intrinsically on the ambient scheme, and it could be defined just alluding to the structure sheaf $\CO_X,$ as it appears for example in \cite{OTT}.The one we expose in the present work is sightly different, but it is handier in our context and it is easy to prove that both definitions are equivalent.
\end{remark}

\begin{paragrafo}{\bf Higher Atiyah classes.}\end{paragrafo}
Given a sheaf $\CL\in\CA_{qc}(X)$ its exterior powers form a graded algebra with multiplication given by the natural maps:
\[\wedge^{i,j}:\bigwedge^i\CL\otimes\bigwedge^j\CL\lto\bigwedge^{i+j}\CL.\]
In particular we have the exterior algebra of the sheaf of differentials. As it is customary, we will denote $\Omega_{X|S}^i:=\bigwedge^i\Omega_{X|S}.$ Note that, by smoothness, all $\Omega_{X|S}^i$ are locally free of finite rank and $0$ if $i>n$. By induction, let $i\geq1$ and assume defined $\at_\CE^i:\CE\to\CE\otimes\Sigma^i\Omega_{X|S}^i.$ Define $\at_\CE^{i+1}$ as the following composition
\[\CE\stackrel{\at_\CE^i}{\lto}\CE\otimes\Sigma^i\Omega_{X|S}^i
\rTo^{\at_\CE^1\ot\id}
\CE\otimes\Sigma^1\Omega_{X|S}^1\otimes\Sigma^i\Omega_{X|S}^i
\stackrel{\wedge^{1,i}}{\lto}
\CE\otimes\Sigma^{i+1}\Omega_{X|S}^{i+1}\]

\begin{lemma}\label{natat}
For every $n\in\mathbb{N}$,
\[\at^i:\id_{\BD(\CA_{qc}(X))}\lto-\otimes\Sigma^i\Omega_{X|S}^i\]
defines a natural transformation.
\end{lemma}
{\bf Proof} Let $q:\CE\to\CE'$ be a morphism of perfect complexes. We have to show that
\[(q\otimes\id_{\Sigma^{i+1}\Omega_{X|S}^{i+1}})\at_\CE^i=\at_{\CE'}^iq.\]
For $i=1$ it holds because the diagram
\begin{diagram}
\CE\otimes\Sigma^1\Omega_{X|S}^1 & \rTo & \CE\otimes\CP_{X|S}^1 & \rTo &
\CE & \rTo^{\at_\CE^1} & \CE\otimes\Sigma^1\Omega_{X|S}^1\\
\dTo^{q\otimes\id} & & \dTo^{q\otimes\id}& & \dTo_q & & \dTo_{q\otimes\id}\\
\CE'\otimes\Sigma^1\Omega_{X|S}^1 & \rTo & \CE'\otimes\CP_{X|S}^1 & \rTo &
\CE' & \rTo^{\at_{\CE'}^1} & \CE'\otimes\Sigma^1\Omega_{X|S}^1
\end{diagram}
commutes. For $i>1$ it follows by induction once we realize that the iteration is clearly built from similarly natural maps.
\qed

As a consequence of the lemma, $\at^i$ constitutes an orientation in the sense of $\S4$ with $C=\Sigma^i\Omega_{X|S}^i$. It gives a character with values in Hodge cohomology, namely
\[\ch_i(\CE):=\Lef(\id_\CE,\at^i):\CO_X\lto\Sigma^i\Omega_{X|S}^i.\]
Observe that $\ch_i(\CE)\in H^i(X,\Omega_{X|S}^i)$.

\begin{proposition}\label{invqischi}
If $\CE$ and $\CE'$ are isomorphic in $\BD(\CA_{qc}(X))$ then
\[\ch_i(\CE)=\ch_i(\CE')\]
for every $i\in\mathbb{N}$.
\end{proposition}

{\bf Proof} We reduce at one to the case in which there is a quasi-isomorphism $q:\CE\to\CE'.$ In this case, it follows from the commutativity of the diagram
\begin{diagram}
                                 &                          &  \CO_X &                    &    \\
                                 &\ldTo^{u_\CE} &       \text{D}_u        &\rdTo^{u_\CE'} &    \\
\CE\otimes \mathrm{D}\CE &\rTo^{q\otimes\id}  &\CE'\otimes\mathrm{D}\CE & \lTo^{\id\otimes\mathrm{D}q}   &\CE'\otimes \mathrm{D}\CE'\\
\dTo^{\gamma}&     {\text{nat}}_\gamma          &\dTo^{\gamma}   &     {\text{nat}}_\gamma &\dTo_{\gamma} \\
\mathrm{D}\CE\otimes\CE &\rTo^{\id\otimes q}  &\mathrm{D}\CE\otimes\CE' & \lTo^{\mathrm{D}q\otimes\id}   &\mathrm{D}\CE'\otimes\CE'\\
\dTo^{\id\otimes\at_{\CE}^i}& {\text{nat}}_{\at} & \dTo^{\id\otimes\at_{\CE'}^i} & \text{triv} &    \dTo_{\id\otimes\at_{\CE'}^i} \\
\mathrm{D}\CE\otimes\CE\otimes\Sigma^i\Omega_X^i &\rTo^{\id\otimes q\otimes\id}  &\mathrm{D}\CE\otimes\CE'\otimes\Sigma^i\Omega_X^i  & \lTo^{\mathrm{D}q\otimes\id\otimes\id}   &\mathrm{D}\CE'\otimes\CE'\otimes\Sigma^i\Omega_X^i\\
\dTo^{\id} & & \dTo^{\mathbb{D}q^{-1}\otimes q^{-1}\otimes\id} & & \dTo_{\id}\\
\mathrm{D}\CE\otimes\CE\otimes\Sigma^i\Omega_X^i &\rTo^{\mathrm{D}q\otimes \id\otimes\id}  &\mathrm{D}\CE'\otimes\CE\otimes\Sigma^i\Omega_X^i  & \lTo^{\id\otimes q\otimes\id}   &\mathrm{D}\CE'\otimes\CE'\otimes\Sigma^i\Omega_X^i\\
                                 &\rdTo_{t_\CE\otimes\id} &   \text{D}_t            &\ldTo_{t_{\CE'}\otimes\id} &
                                 \\
                                &                          &    \Sigma^i\Omega_X^i&                    &
\end{diagram}

Diagrams labeled by $\text{D}_u$ and $\text{D}_t$ commute by Lemma \ref{diagdual}. Those labeled by $\text{nat}_{\gamma}$ or $\text{nat}_{\at}$ by the naturality of $\gamma$ and $\at^i$ respectively. The subdiagram labeled by triv commutes by obvious reasons. The remaining squares commute by an obvious computation keeping in mind that $(\mathrm{D}q)^{-1}=\mathrm{D}q^{-1}$
\qed

{\it Remark}. We recall that even if we are able to represent $a$ by a morphism of complexes this not need to be the case for $q^{-1}.$

\begin{theorem}\label{addchi} Let
\[\CE\stackrel{f}{\lto}\CF\stackrel{g}{\lto}\CG\stackrel{h}{\lto}\Sigma\CE\]
be a distinguished triangle. Then
\[\ch_i(\CF)=\ch_i(\CE)+\ch_i(\CG)\]
for every $i\in\mathbb{N}$.
\end{theorem}

{\bf Proof} In virtue of Proposition \ref{invqischi} we can assume that $\CG=\cone(f).$ By definition $\ch_i(-)$ it must be checked the commutativity of the outer maps in a
 diagram analogous to the one considered in the proof of Theorem \ref{addtraces} and, as we done there, the strategy consists on dividing it in several smaller subdiagrams. Namely we consider
\begin{diagram}
                                 &                          &  \CO_X &                    &    \\
                                 &\ldTo^{u_\CG,u_E} &               &\rdTo^{u_F} &    \\
(\CG\otimes \mathrm{D}\CG)\oplus(\CE\otimes \mathrm{D}\CE) &  & &    &(\CF\otimes \mathrm{D}\CF)\\
\dTo^{(\gamma,\gamma)}&                              &   &     &\dTo_{\gamma} \\
(\mathrm{D}\CG\otimes\CG)\oplus(\mathrm{D}\CE\otimes\CE) &  & &    &(\mathrm{D}\CF\otimes\CF) \\
\dTo^{(\id\otimes\at^i(\CG),\id\otimes\at^i(\CE))}&                              &   & &    \dTo_{\id\otimes\at^i(\CF)} \\
(\mathrm{D}\CG\otimes\CG\otimes \Sigma^i\Omega_X^i)\oplus(\mathrm{D}\CE\otimes\CE\otimes \Sigma^i\Omega_X^i) &  & &    &\mathrm{D}\CF\otimes\CF\otimes \Sigma^i\Omega_X^i \\
                                 &\rdTo_{(t_\CG,t_\CE)} &               &\ldTo_{t_\CF} &    \\
                                &                          &   \Sigma^i\Omega_X^i&                    &
                                \end{diagram}

Taking the an object $\mathcal{W}$ in the fashion of (\ref{tc3h}), as the Atiyah class is a natural transformation, we can proceed as in the proof of Theorem \ref{addtraces}, choosing in the present case the morphism $m:=\at^i(\mathcal{W}):\mathcal{W}\to\mathcal{W}\otimes\Sigma^i\Omega_X^i$ to conclude.

\begin{remark}
As we have seen, for perfect objects, the Atiyah class depends only on the scheme $X,$ in the sense that we just tensorize by the concrete object, but there is no braid with objects arising in the definition of the Atiyah class of $X$, so there is a commutative diagram
\begin{diagram}
\CE & \rTo^{f} & \CF \\
\dTo^{\at(\CE)} & & \dTo_{\at(\CF)}\\
\CE\otimes\oplus_n(\Sigma^n\Omega_X^n) & \rTo_{f\otimes\id} & \CF\otimes\oplus_n(\Sigma^n\Omega_X^n)
\end{diagram}
Then the Atiyah class is functorial.
\end{remark}

\begin{paragrafo}{\bf The Atiyah orientation and the Chern character.}
\end{paragrafo}
Suposse further that $n!$ is a unit in the ring of sections of $\CO_S$ and, as a consequence, also in $\CO_X.$ For instance it is enough that there exists morphism $S\to \spec \mathbb{F}_p$ with $p$ a prime integer greater than $n$ or a morphism $S\to \spec \mathbb{Q}.$ If $\CE\in\BD(\CA_{qc}(X))$ is a perfect complex, we define its total Atiyah class
\[\At_\CE:\CE\lto\CE\otimes\bigoplus_{i=0}^n\Sigma^i\Omega_{X|S}^i\]
as
\[\At_\CE:=\sum_{i=0}^n\frac{1}{i!}\at_\CE^i\]
It follows from Lemma \ref{natat} that
\[\At:\id_{\BD(\CA_{qc}(X))}\lto-\otimes\bigoplus_{i=0}^n\Sigma^i\Omega_{X|S}^i\]
is a natural transformation, and therefore constitutes an orientation in the sense of $\S$4. In this case $C=\bigoplus_{i=0}^n\Sigma^i\Omega_{X|S}^i$, so we can define the Chern character with values in Hodge cohomology as
\[\ch(\CE):=\Lef(\id_\CE,\At):\CO_X\lto\bigoplus_{i=0}^n\Sigma^i\Omega_{X|S}^i.\]
Note that $\ch(\CE)=\sum_{i=0}^n\frac{1}{i!}\ch_i(\CE)$ and that $\ch(\CE)\in\bigoplus_{i=0}^n H^i(X,\Sigma^i\Omega_{X|S}^i).$

\begin{corollary}\label{invqisch}
Suppose that $n!$ is a unit in the ring of sections of $\CO_S$. If $\CE$ and $\CE'$ are isomorphic in $\BD(\CA_{qc}(X))$ then
\[\ch(\CE)=\ch(\CE')\]
\end{corollary}

{\bf Proof}
It is an inmediate consequence of Proposition \ref{invqischi} and the definition or $\ch$.
\qed

\begin{corollary}\label{addch}
Suppose that $n!$ is a unit in the ring of sections of $\CO_S$ and let
\[\CE\stackrel{f}{\lto}\CF\stackrel{g}{\lto}\CG\stackrel{h}{\lto}\Sigma\CE\]
be a distinguished triangle. Then
\[\ch(\CF)=\ch(\CE)+\ch(\CG).\]
\end{corollary}
{\bf Proof}
Again, an  inmediate consequence of Theorem \ref{addchi}.
\qed

The previous corollaries may be organized as a map that plays a basic role in Riemann-Roch theorems.

Denote by $\BD(X)_{cp}$ the full subcategory of $\BD(\CA_{qc}(X))$ whose objects are perfect complexes. This category is skeletally small (see \cite[Appendix F]{TT}), so there exist  a set $U$ containing at least a representative for every class of isomorphism of objects in $\BD(X)_{cp}.$ We consider $K^0(X)$ the quotient of the free group over $U$ by the subgroup generated by the relations $\CE-\CE'$ if $\CE\cong\CE'$ and the Euler equation $\CE-\CF-\CG$ for every distinguished triangle $\CE\to\CF\to\CG\stackrel{+}{\to}$ with vertices in $U$ (see \cite{SGA6}) for details).
Whit this notation we have the following

\begin{theorem}In the previous setting, suppose that $n!$ is a unit in the ring os sections of $\CO_S$. There is a homomorphism of groups
\[\ch:K^0(X)\lto\bigoplus_{i=0}^n H^i(X,\Omega_{X|S}^i).\]
\end{theorem}
{\bf Proof}
It is clear that the Chern character $\ch$ defines a map with source $U$ that may be extended to the free group it generates by linearity. Corollary \ref{invqisch} guarantees that the relations for isomorphisms go to $0$ in the target and Corollary \ref{addch} does the same for the Euler relations. Whence, the conclusion follows.
\qed

{\bf Remark.} This map plays the same role as the Chern character in intersection theory, with the $i^{\text{th}}$ Chow group replaced by the $i^{\text{th}}$ Hodge cohomology group $H^i(X,\Omega_{X|S}^i)$.


\begin{thebibliography}{ABCDD}

\bibitem[Be]{Be} Beke, T. :{\it Sheafifiable homotopy model categories}, Math. Proc. Camb. Phil. Soc. {\bf 129}, (2000), no. 3, p. 447-475.

\bibitem[SGA6]{SGA6} Berthelot, P.; Grothendieck, A.; Illusie, L.:{\it Th\'eorie des Intersections et Th\'eor\'eme de Riemann-Roch (SGA 6)}, Lecture Notes in Math., Vol. {\bf 225}, Sp`ringer-Verlag, New York, (1971).

\bibitem[CD]{CD} Cisinski, D. C.; Deglise, F.:{\it Local and stable homological algebra in  Grothendieck abelian categories}, Homology, Homotopy Appl. {\bf 11} (2008), no. 1, 219-260.

\bibitem[EGA IV]{EGA IV} Dieudonn\'e, J; Grothendieck, A.; Illusie, L.:{\it Fundaments de geom\'etrie alg\'ebrique IV. Et\'ude locale des schemas et des morphismes de schemas, Quatri\'eme partie}, Publications mathematiques de l'I.H.ƒ.S., {\bf 32}, (1967), 5-361.


\bibitem[EK]{Ei} Eilenberg, S.; Kelly, G. M.: Closed categories. Proc. Conf. Categorical Algebra
    (La Jolla, Calif., 1965), Springer, New York, (1966), 421-562.


\bibitem[F]{Fe} Ferrand, D.:{\it On the additivity of the trace in derived categories}, (2005). Available at arXiv:math.CT/0506589.

\bibitem[FM]{FM} Fulton, W.; MacPhenson, R.:{\it Categorical framework for the study of singular spaces}, Mem. Amer. Math. Soc. {\bf 31}, (1981), no. 243.

\bibitem[HPS]{HPS} Hovey, M.; Palmieri, J. H.; Strickland, N. P.: {\it Axiomatic stable homotopy theory}. Mem. Amer. Math. Soc. {\bf 128} (1997), no. 610.

\bibitem[KN]{Kell-Nee} Keller, B., Neeman, A.:{\it The connection between May's Axioms for a triangulated tensor product and Happel's description of the derived category of the quiver $D_4.$}, Documenta Mathematica .{\bf 7}, (2002), 535-560.

\bibitem[LMSMc]{LMS} Lewis, L. G.,Jr.; May, J. P.; Steinberger, M.: McClure, J. E.:{\it Equivariant sable homotopy theory}. Lecture Notes in Math. {\bf 1213}. Springer-Verlag, Berlin, (1986).

\bibitem[Lip]{Lip} Lipman, J.: {\it Notes on derived categories and Grothendieck Duality for Diagrams of Schemes}, Lecture Notes in Math., {\bf 1960}, Springer-Verlag, (2009).

\bibitem[M]{May} May, J.P.: {\it The Additivity of Traces in Triangulated Categories}, Advances in  Mathematics {\bf 163}, (2001), 34-73.

\bibitem[N]{Nee2} Neeman, A.: \textit{Triangulated Categories}, Annals of mathematics studies, {\bf 148}, Princeton Universitity Press, (2001).

\bibitem[N2]{Neeperf} Neeman, A.: {\it Derived categories and Grothendieck duality}, CRM Preprints {\bf 791}, (2008), avaliable at

    http://crm.es/Publications/08/Pr791.pdf.


\bibitem[OTT]{OTT} O'Brian, N. R.; Toledo, D.;Tong, Y.L.: {\it The trace map and characteristic classes for coherent sheaves}. Amer. J. Math. {\bf 103} (1981), no. 2, 225-252.

\bibitem[Se]{Se} Serp\'e, C.: {\it Resolution of unbounded complexes in Grothendieck categories}. J. of Pure and Applied Algebra {\bf 177}, (2003), 102-112.


\bibitem[Sp]{Sp} Spalstestein, N.: {\it Resolutions of unbounded complexes}. Compositio Math. {\bf 65}, (1988), no. 2, 121-154.


\bibitem[TT]{TT}Thomason, R. W.; Trobaugh, T.: {\it Higher algebraic K-theory of schemes and of derived categories. The Grothendieck Festschrift, Vol. III.} Prog.Math., {\bf 88}, (1990), 247-435.

\end{thebibliography}
\end{document}